\documentclass[12pt]{amsproc}
\usepackage{amssymb}
\usepackage{epsf}

\newfam\cyrfam
\font\Cyr=wncyr10 at 8pt   
\font\CYR=wncyr10 at 12pt  

\copyrightinfo{2000}{(Frank Sottile)}

\newtheorem{theorem}{Theorem}[section]
\newtheorem{lemma}[theorem]{Lemma}
\newtheorem{prop}[theorem]{Proposition}

\theoremstyle{definition}

\theoremstyle{remark}
\newtheorem{remark}[theorem]{Remark}

\theoremstyle{remark}
\newtheorem{remarks}[theorem]{Remarks}

\numberwithin{equation}{section}

\begin{document}

\title[Rational Curves on Grassmannians]%
      {Rational Curves on Grassmannians:\\
       systems theory, reality, and transversality}

\author{Frank Sottile}
\address{Department of Mathematics and Statistics\\
        University of Massachusetts\\
        Amherst, Massachusetts 01003\\
        USA}
\email[Frank Sottile]{sottile@math.umass.edu}
\urladdr[Frank Sottile]{http://www.math.umass.edu/\~{}sottile}
\thanks{Research supported in part by NSF grant DMS-0070494}
\thanks{Based upon a talk by the author in the Special Session on Enumerative
Geometry in Physics at the AMS sectional meeting in Lowell, Massachusetts,
April 1-2, 2000.}


\begin{abstract}
We discuss a particular problem of enumerating rational curves on a
Grassmannian from several perspectives, including systems theory, real
enumerative geometry, and symbolic computation.
We also present a new transversality result, showing this problem is enumerative
in all characteristics.

While it is well-known how this enumerative problem arose in mathematical
physics and also its importance to the development of quantum cohomology, it is  
less known how it arose independently in mathematical systems theory.
We describe this second story.
\end{abstract}

\maketitle

\section{Introduction}

The enumerative geometry of curves on algebraic varieties has become an
important theme in algebraic geometry.
One motivation for this development was to understand (and prove) remarkable
formulae from theoretical
physics, including a formula of Vafa and Intriligator~\cite{In91,Va92}
involving curves on Grassmannians. 
The story of this direct influence of theoretical physics on algebraic
geometry is well-known.
What is less known is how the problem of enumerating
rational curves on Grassmannians also arose and was solved in systems theory.
Our purpose is to make that story better known and to relate the
different solutions, from physics and from systems theory, of this enumerative problem.
We also discuss some related work in algebraic geometry inspired by systems theory.

We describe this enumerative problem.
Let $m,p\geq 1$ be integers.
The space ${\mathcal M}^q_{m,p}$ of maps $M$ of degree $q$  
from ${\mathbb P}^1$ to ${\it Grass}(p,{\mathbb C}^{m+p})$, 
the Grassmannian of $p$-planes in ${\mathbb C}^{m+p}$, has dimension 
$N:=q(m+p)+mp$~\cite{Clark,St87}.
Given a point $s\in{\mathbb P}^1$ and an $m$-plane $L$ in 
${\mathbb C}^{m+p}$, the set of maps $M$ which satisfy 
$M(s)\cap L\neq\{0\}$ (the $p$-plane $M(s)$ meets the $m$-plane $L$
non-trivially) is a divisor on this space of maps.
We consider the following enumerative problem:
\medskip

\noindent{\bf Question 1.}
{\it 
  Given general points $s_1,s_2,\ldots,s_N\in{\mathbb P}^1$ and general 
  $m$-planes 
  $L_1$, $L_2,\ldots$, $L_N\subset{\mathbb C}^{m+p}$, how many degree $q$ maps 
  $M\colon{\mathbb P}^1\to{\it Grass}(p,{\mathbb C}^{m+p})$ satisfy
  \begin{equation}\label{condition}
   M(s_i)\cap L_i\ \neq\ \{0\}\qquad\mbox{for }\ i=1,2,\ldots, N\ ?
  \end{equation}}

This is a special case of the more general enumerative problem
considered by Vafa and Intriligator~\cite{In91,Va92} who replaced the 
Schubert condition $M(s)\cap L\neq\{0\}$ by general Schubert
conditions and the map 
$M\colon({\mathbb P^1},s_1,s_2,\ldots,s_N)\to{\it Grass}(p,{\mathbb C}^{m+p})$
by a map of a general pointed curve.
There, a formula was proposed involving residues.
This formula was justified by Siebert and Tian~\cite{ST97} by 
computing the (small) quantum cohomology ring of the
Grassmannian, whose structure was also conjectured by Vafa and 
Intriligator.
We describe this part of our story in Section~\ref{sec:last}.

A completely different approach (and motivation) to this enumerative problem
came from systems theory.
Briefly, conditions of the form $M(s)\cap L\neq\{0\}$
arise in the problem of stabilizing a given linear system using dynamic output
compensation~\cite{Ro94}.
In the critical dimension when there are finitely many compensators, the
problem of enumeration was solved by Ravi, Rosenthal, and
Wang~\cite{RRW96,RRW98}, 
who gave the closed formula for the intersection number
$d(m,p;q)$ of Question 1:
 \begin{equation}\label{formula}
  (-1)^{q(p{+}1)}\;N!\ \cdot\
  \sum_{\nu_1+\cdots+\nu_p=q}
  \frac{\prod_{j<k}(k{-}j{+}(\nu_k{-}\nu_j)(m{+}p))}
  {\prod_{j=1}^p(m{+}j{+}\nu_j(m{+}p){-}1)!}\ .
 \end{equation}
One of their motivations was to determine when this number is odd, for then 
there exists 
a real compensator stabilizing a given real linear system.
We describe how this problem in systems theory is a special case of the 
general enumerative problem described above, and also how Ravi, Rosenthal,
and Wang solved this enumeration in Section~\ref{sec:dynamic}. 
We remark that the quantum cohomology of the Grassmannian also has applications
to matrix interpolation problems~\cite{BR93,Ra99}.

The geometric formulation from systems theory (and ideas from numerical
homotopy continuation~\cite{AG97}) were exploited to prove the
following result in real enumerative geometry:
There exist real points $s_1,s_2,\ldots,s_N\in{\mathbb P}^1_{\mathbb R}$ and real
$m$-planes $L_1,L_2,\ldots,L_N\subset{\mathbb R}^{m+p}$ such that there are
$d(m,p;q)$ rational maps 
$M\colon{\mathbb P}^1\to{\it Grass}(p,{\mathbb C}^{m+p})$ of degree $q$
satisfying~(\ref{condition}), and each of these maps is real~\cite{So00}.
Thus the enumerative problem of Question 1 is fully real (in the sense
of~\cite{So97c}). 
A variant of this argument gives the new result that 
Question~1 makes enumerative sense in any characteristic: 
If ${\mathbb K}$ is any algebraically closed field, then for general 
points $s_1,s_2,\ldots,s_N\in{\mathbb P}^1_{\mathbb K}$ and general
$m$-planes $L_1,L_2,\ldots,L_N\subset{\mathbb K}^{m+p}$ there are 
exactly $d(m,p;q)$ degree $q$ rational maps 
$M\colon{\mathbb P}^1\to{\it Grass}(p,{\mathbb K}^{m+p})$ 
satisfying~(\ref{condition})~\cite{So_trans}.
The point here is that the corresponding varieties intersect transversally and
so the solutions occur without multiplicities.
We give a proof of these results in Section~\ref{sec:reality}, where we also solve the
enumerative problem of Question 1 without reference to the Chow or quantum Chow
rings, the usual tools of enumerative geometry. 

Ravi, Rosenthal, and Wang~\cite{RRW96,RRW98} also showed that $d(m,p;q)$
equals the number of saturated chains in a certain poset of quantum Pl\"ucker
coordinates.
This is the degree of the singular Uhlenbeck
compactification~\cite{SU81,BDW96} of the 
space of rational curves in the Grassmannian in a natural projective
embedding, also called the quantum Grassmannian.
Its degree may be computed from its defining ideal.
In~\cite{SS_SAGBI}, quantum Pl\"ucker relations for this ideal were
constructed, giving a different proof that this degree equals the number of
chains in the poset of quantum Pl\"ucker coordinates.
We describe that in Section~\ref{sec:eqs} and give another proof that $d(m,p;q)$ equals the
number of chains in that poset.

In the last section, we not only describe some of the classical story motivated by
physicists, but also relate the formula~(\ref{formula}) of Ravi, Rosenthal, and Wang to
the formula of Vafa and Intriligator.
This involves another, intermediate formula~(\ref{fnD}).
We conclude by discussing some further aspects of the quantum cohomology ring of the
Grassmannian,  including how it arose in representation theory and open problems involving
quantum Littlewood-Richardson numbers.


\section{Dynamic Control of Linear Systems}\label{sec:dynamic}
In control theory, the greatest interest is to obtain results valid over the
real numbers ${\mathbb R}$.
As in algebraic geometry, the strongest and most elegant results are true
only for the complex numbers ${\mathbb C}$.
Also as in algebraic geometry, much of the theory may be developed over any
field.
To that end, we let ${\mathbb K}$ denote an arbitrary field, keeping in mind
the special cases of when ${\mathbb K}={\mathbb R}$ or 
${\mathbb K}={\mathbb C}$.\smallskip

Suppose we have a time-invariant physical system with $m$ inputs 
${\bf u}\in{\mathbb K}^m$ 
and $p$ outputs ${\bf y}\in{\mathbb K}^p$ whose evolution is
governed by a system of constant coefficient linear differential equations
$$
   0\ =\ F({\bf u},{\bf u}',\ldots;{\bf y},{\bf y}',\ldots)\,.
$$
One important way in which such a linear system arises is from a linear
perturbation of a non-linear system.

Introducing auxiliary variables or {\it internal states}
${\bf x}\in{\mathbb K}^n$, we can transform this into a first order system of
linear differential equations
 \begin{equation}\label{system}
  \begin{array}{rcl}
   \frac{d}{dt}{\bf x}&=&A{\bf x}+B{\bf u}\\
               {\bf y}&=&C{\bf x}+D{\bf u}\,,
  \end{array}
 \end{equation}
where $A,B,C$, and $D$ are matrices of the appropriate size.
The matrix $D$ represents a direct linear dependence of ${\bf y}$ on 
${\bf u}$.
Systems with $D=0$, where the dependence of ${\bf y}$ on ${\bf u}$
is purely dynamic, are called {\it strictly proper}.
The representation~(\ref{system}) is called a {\it state space form} or
{\it state space realization} of the original system. 
There are many ways to realize a given system in state-space form and a
fundamental invariant, the {\it McMillan degree}, is the minimal number $n$ of
internal states needed to obtain such a first order linear evolution
equation.
The McMillan degree measures the complexity of a linear system.

A system is {\it observable} if the joint kernel of the matrices $CA^k$
for $0\leq k< n$ is zero, which implies that the internal states (${\bf x}$)
may be recovered from knowledge of ${\bf y}(t)$ and ${\bf u}(t)$.
It is {\it controllable} if the matrices 
$A^kB$ for $0\leq k<n$ span ${\mathbb K}^n$, which implies that the system
may be driven to any fixed internal state.
A state space realization~(\ref{system}) of a system is {\it minimal} ($n$ is
its McMillan degree) if and only if it is both observable and
controllable~\cite[\S13]{Fa99}.

\subsection{Rational curves on Grassmannians}
We give another fundamental representation of a linear system that links
systems theory to the (quantum) cohomology of the Grassmannian.
Consider the Laplace transform of~(\ref{system})
$$
  \begin{array}{rcl}
    s\cdot {\bf x}(s)&=&A\cdot{\bf x}(s)+B\cdot{\bf u}(s)\,,\\
           {\bf y}(s)&=&C\cdot{\bf x}(s)+D\cdot{\bf u}(s)\,.
  \end{array}
$$
We eliminate ${\bf x}$ and solve
$$
  {\bf y}(s)\ =\ \left(C(sI_n-A)^{-1}B + D\right){\bf u}(s)\,.
$$
This $p$ by $m$ matrix $\Gamma(s):=C(sI_n-A)^{-1}B + D$ of rational functions
is called the {\it transfer function} of the original system.
It represents the response of the system in the frequency domain.

The transfer function determines a curve in 
${\it Grass}(p,{\mathbb K}^{m+p})$ by
$$
   {\mathbb P}^1\ni s\ \longmapsto\ \mbox{\rm column space}
    \left[\begin{array}{c}I_m\\\Gamma(s)\end{array}\right]\,,
$$
whenever this is well-defined.
This {\it Hermann-Martin} curve extends to ${\mathbb P}^1$ and its degree 
is equal to the McMillan degree of the system.
Recall that the degree of a curve 
 $M\colon{\mathbb P}^1\to{\it Grass}(p,{\mathbb K}^{m+p})$
 has three equivalent descriptions:
\begin{enumerate}
\item The number of points $s\in{\mathbb P}^1$ such that 
      $M(s)\cap L\neq\{0\}$, where $L$ is a general $m$-plane.
\item The maximum degree of the (rational-function) minors of
      any $(m+p)$ by $p$ matrix of rational functions whose column space
      gives the map $M$.
\item The degree of the pullback of the generator ${\mathcal O}(1)$ of
      the Picard group of ${\it Grass}(p,{\mathbb K}^{m+p})$.
\end{enumerate}

One concrete way to see that the transfer function defines a curve in the
Grassmannian is via the algebra of polynomial matrices.
A matrix $\Gamma(s)$ of rational functions is {\it proper} if 
$\lim_{s\to \infty}\Gamma(s)$ exists and {\it strictly proper} if that limit is
zero.
The transfer function of the linear system~(\ref{system}) is proper,
since $\lim_{s\to\infty}\Gamma(s)=D$, and strictly proper linear
systems have strictly proper transfer functions.
Given a proper matrix of rational functions $\Gamma(s)$ of size $p$ by $m$,
consider factorizations
$$
   \Gamma(s)\ = \ P(s)Q(s)^{-1}
$$
where $P(s)$ is a $p$ by $m$ matrix of polynomials and $Q(s)$ is a 
$m$ by $m$ matrix of polynomials with non-zero determinant.
There are many ways to do this:
One could, for instance, let $Q(s)$ be the diagonal matrix with entries
$f(s)$, the least common multiple of the denominators of the entries of
$\Gamma(s)$.
There is a unique minimal, or (right) coprime factorization.

\begin{theorem}\label{thm11}
 Suppose $\Gamma(s)$ is a proper transfer function of a linear
 system of McMillan degree $n$.
 Then there exist matrices $P(s),Q(s)$ of polynomials such that 
 \begin{enumerate}
  \item[({\it i})]
        $P(s)$ and $Q(s)$ are coprime in that there exist matrices of 
        polynomials $X(s)$ and $Y(s)$ satisfying
$$
       X(s) Q(s) +  Y(s) P(s)\ =\ I_m\,.
$$

  \item[({\it ii})]
         Any other factorization $\Gamma(s)=N(s)D(s)^{-1}$ into matrices of
        polynomials has
$$
   \deg \det D(s)\  \geq\ \deg\det Q(s)\ =\ n\,.
$$

  \item[({\it iii})]
         $P(s)$ and $Q(s)$ are unique up to multiplication on the right by
        elements of $GL_m({\mathbb K}[s])$.
\end{enumerate}
\end{theorem}

Theorem~\ref{thm11} is proven, for instance in any of~\cite{Ro70} 
or~\cite[\S22]{De88} or~\cite[\S4]{Fa99}.  

By ({\it i}) and the factorization, 
$P(s)Q(s)^{-1}=C(sI_n{-}A)^{-1}B+D$, the determinants of $Q(s)$ and of
$sI_n{-}A$ have the same roots.
We call $P(s)Q(s)^{-1}$ a {\it right coprime factorization}
of $\Gamma(s)$.
By ({\it i}), the Hermann-Martin curve is also represented by
 \begin{equation}\label{H-M}
   s\ \longmapsto\ \mbox{\rm column space}\left[
        \begin{array}{c}Q(s)\\P(s)\end{array}\right] \,,
 \end{equation}
which has dimension $m$ for all $s\in{\mathbb P}^1$, as $\Gamma(s)$ is proper.
Since $\lim_{s\to\infty}\Gamma(s)=D$, the value of the curve at infinity is
the column space of 
$$
  \lim_{s\to\infty} 
    \left[\begin{array}{c}Q(s)\\P(s)\end{array}\right]
  \quad\ =\ \quad
  \left[ \begin{array}{c}I_m\\D\end{array}\right]\,.
$$
Thus the maximal minors of the matrix
$$
\left[\begin{array}{c}Q(s)\\P(s)\end{array}\right]
$$
have degree at most the degree of the principal minor
$\det Q(s)$, which is $n$.
This shows that the Hermann-Martin curve has degree $n$.
In this way, a linear system~(\ref{system}) with $m$ inputs and $p$ outputs
of McMillan degree $n$ corresponds  to a rational curve
$M\colon{\mathbb P}^1\to\mbox{\it Grass}(p,{\mathbb K}^{m+p})$ of degree $n$.
In fact, every such rational curve comes from a linear system~\cite{HM78}.

An informal way to see this is to first observe that the entries of the
matrices 
$A,B,C$, and $D$ in~(\ref{system}) give the set of all possible state-space
realizations of $m$-input $p$-output linear systems with $n$ internal states
the structure of affine space of dimension $n^2+nm+np+mp$.
The conditions of controllability and observability for the system to be
minimal are the non-vanishing of certain polynomials in the entries of
$A,B,C$, and so the set of all such systems of McMillan degree $n$ is an
open subset of this affine space.
Changing coordinates of the internal states ${\bf x}$ gives a free 
$GL_n({\mathbb K})$-action on these minimal realizations whose orbits are
exactly the fibres of the map
$$
  \left\{
     \mbox{\rm Minimal state-space realizations}
  \right\}
     \quad\longrightarrow\quad
  \left\{
     \mbox{\rm Proper transfer functions}
  \right\} \ .
$$
Thus the space of Hermann-Martin curves of $m$-input $p$-output linear
systems of McMillan degree $n$ has dimension $nm+np+mp$,
which is equal to the dimension of the space ${\mathcal M}^q_{p,m}$ 
of degree $n$ rational maps 
to $\mbox{\it Grass}(m,{\mathbb K}^{m+p})$.

In fact, the Hermann-Martin curves constitute an open subset of this space
of rational curves, and there are
very natural objects from systems theory that yield the full space of rational
curves, as well as various compactifications of this space.
The work of Hermann and Martin~\cite{HM78} continued work of Clark~\cite{Clark},
who showed that the space of transfer functions is a smooth manifold.
Later, Helmke~\cite{He86} studied topological properties of this space and
Hazewinkel~\cite{Ha80} and Byrnes~\cite{By83} studied compactifications of this
space.
This work was revived by Rosenthal, who introduced the quantum Grassmannian into
systems theory in his 1990 PhD thesis~\cite{Rosen_PhD}.
See~\cite{RR94} for a discussion and further references.


\subsection{Feedback control and Schubert calculus}
Given a strictly proper linear system
 \begin{equation}\label{spls}
  \begin{array}{rcl}
   \frac{d}{dt}{\bf x}&=&A{\bf x}+B{\bf u}\,,\\
               {\bf y}&=&C{\bf x}\,,
  \end{array}
 \end{equation}
we would like to control its behavior
using dynamic output feedback.
That is, we couple its inputs ${\bf u}$ to its outputs ${\bf y}$ through a $p$-input,
$m$-output linear system of McMillan degree $q$, called a {\it dynamic compensator}.
Consider a minimal state-space realization of this compensator
 \begin{equation}\label{compensator}
  \begin{array}{rcl}
   \frac{d}{dt}{\bf z}&=&F{\bf z}+G{\bf y}\\
               {\bf u}&=&H{\bf z}+K{\bf y}\,,
  \end{array}
 \end{equation}
where ${\bf z}\in{\mathbb K}^q$ are the internal states, and $F,G,H$, and
$K$ are matrices of the appropriate size, $K$ representing a constant
(residual) linear feedback law.

Schematically we have:
$$
\begin{picture}(363,100)
\thicklines

\put(0,70){Given system:}
\put(0,20){Compensator:}

\put(100,73){\line(1,0){100}}
   \put(200,73){\line(-1,1){10}}  \put(200,73){\line(-1,-1){10}}
\put(135,78){${\bf u}\in{\mathbb K}^m$}
\put(200,60){\line(0,1){25}} \put(253,60){\line(0,1){25}}
\put(200,60){\line(1,0){53}} \put(200,85){\line(1,0){53}} 
\put(210,70){${\bf x}\in{\mathbb K}^n$}
\put(253,73){\line(1,0){100}}
   \put(353,73){\line(-1,1){10}}  \put(353,73){\line(-1,-1){10}}
\put(288,78){${\bf y}\in{\mathbb K}^p$}

\put(150,73){\line(0,-1){50}}   \put(200,23){\line(-1,0){50}}
   \put(150,73){\line(1,-1){10}}  \put(150,73){\line(-1,-1){10}}
\put(200,10){\line(0,1){25}} \put(253,10){\line(0,1){25}}
\put(200,10){\line(1,0){53}} \put(200,35){\line(1,0){53}} 
\put(210,20){${\bf z}\in{\mathbb K}^q$}
\put(303,73){\line(0,-1){50}}  \put(253,23){\line(1,0){50}}
   \put(253,23){\line(1,1){10}}  \put(253,23){\line(1,-1){10}}

\end{picture}
$$

We obtain a closed-loop or autonomous system from~(\ref{system})
and~(\ref{compensator}) by eliminating ${\bf y}$ and ${\bf u}$
\begin{equation}\label{autonomous}
  \frac{d}{dt}\left[\begin{array}{c}{\bf x}\\{\bf z}\end{array}\right]
   \ = \ 
  \left[\begin{array}{cc}A+BKC&BH\\GC&F\end{array}\right]\ 
  \left[\begin{array}{c}{\bf x}\\{\bf z}\end{array}\right]\ .
\end{equation}

The behavior of this autonomous system is determined by the $n+q$
eigenvalues of the matrix, that is, by the zeroes of the (monic)
characteristic polynomial
 \begin{equation}\label{charpoly}
   \varphi(s)\ :=\ \det\left( s I_{n+q} -
   \left[\begin{array}{cc}A+BKC&BH\\GC&F\end{array}\right]\right)\ .
 \end{equation}

The pole placement problem asks the inverse question:
\medskip
 
\noindent{\bf Pole Placement Problem.}
 Given a strictly proper $m$-input $p$-output linear system of McMillan
 degree $n$~(\ref{spls}) and a desired behavior represented by a monic
 characteristic 
 polynomial $\varphi(s)$ of degree $n+q$, for which dynamic
 compensators~(\ref{compensator}) does the corresponding autonomous
 system~(\ref{autonomous}) have characteristic polynomial $\varphi$?
\medskip

The reason for the word pole is that the zeroes of the characteristic polynomial
are the poles of a transfer function.
A linear system of McMillan degree $n$ is {\it arbitrarily pole-assignable
by degree $q$ compensators} (over ${\mathbb K}$) 
if the pole placement problem may be solved for all monic polynomials
$\varphi$ of degree $n+q$.

\begin{remark}
Pole placement is a fundamental design problem for linear systems.
When ${\mathbb K}={\mathbb R}$, an 
important property of an autonomous real linear system is whether or not it
is stable, that is, whether or not all of the roots of its characteristic
polynomial have negative real parts.
In other situations, the control engineer may wish to destabilize a system.
For discrete-time systems (which have an identical formalism), stability is
achieved by placing the roots of the characteristic polynomial on the unit
circle.
These questions of placing poles in subsets of the complex plane are strictly
weaker than the pole placement problem, yet little is known about them.
Here is an important related question concerning stability.
\medskip

\noindent{\bf Minimal Stability.}
 Given a strictly proper $m$-input $p$-output real linear system of McMillan
 degree $n$, what is the minimal McMillan degree $q$ of a real dynamic
 compensator~(\ref{compensator}) for which the corresponding autonomous 
 system~(\ref{autonomous}) is stable?\medskip
\end{remark}

When ${\mathbb K}$ is algebraically closed, 
the pole placement problem may be solved for $q\geq n-1$~\cite{BP70} and 
$q\leq (n-mp)/(m+p-1)$ is necessary~\cite{Ro94} and sufficient~\cite{WH78} for
generic systems. 
Thus for $q$ large enough there exist stabilizing dynamic
compensators.
The minimal stability problem is particularly important when the original 
system arises as a linear perturbation of a non-linear system.
In this case, it asks how cheaply may we damp linear perturbations of
McMillan degree $n$.
\medskip

We investigate the pole placement problem.
Given a strictly proper system~(\ref{spls}) and a monic characteristic polynomial 
$\varphi(s)$ with distinct roots $s_1,s_2\ldots,s_{n+q}$, we seek matrices
$F,G,H$, and $K$ for which 
$$
  \det\left( s_i I_{n+q} -
  \left[\begin{array}{cc}A+BKC&BH\\GC&F\end{array}\right]\right)\ =\ 0
  \qquad\mbox{for}\quad i=1,2,\ldots,n+q\,.
$$
This gives $n+q$ equations in the $q^2+pq+mq+mp$ entries of 
$F,G,H$, and $K$.  
Since $GL_q({\mathbb K})$ acts on these data, giving equivalent systems and
fixing $\varphi$, we expect that the pole placement problem is solvable over
the complex numbers when  
$$
   n+q \ \leq\ pq + mq + mp\,.
$$
This is in fact the case for generic systems, as we shall see.

We reformulate the dynamic pole placement problem geometrically.
Each step below involves only row or column operations applied to the
matrix involved.
\begin{eqnarray*}
  \varphi(s)&=& 
   \det\left[\begin{array}{cc}sI_n-A-BKC&-BH\\-GC&sI_q-F\end{array}\right]\\
  &=&
   \det\left[\begin{array}{cccc}
        sI_n-A-BKC&  -BH & BK& -B \\
           -GC    &sI_q-F& 0 &  0 \\
            0     &   0  &I_p&  0 \\
            0     &   0  & 0 & I_m  \end{array}\right]  \\
  &=&
   \det\left[\begin{array}{cccc}
        sI_n-A &   0  & 0 & -B \\
           0   &sI_q-F& G &  0 \\
           C   &   0  &I_p&  0 \\
           0   &  -H  & K & I_m  \end{array}\right]  \\
  &=&
   \begin{array}{l} \det[sI_n-A]\\
                    \times \det[sI_q-F] \end{array}
    \times 
    \det\left[\begin{array}{cccc}
        I_n & 0 &       0      &-(sI_n-A)^{-1}B \\
         0  &I_q&(sI_q-F)^{-1}G&       0        \\
         C  & 0 &      I_p     &       0        \\
         0  &-H &       K      &      I_m    \end{array}\right]  
\end{eqnarray*}
\noindent 
This becomes
$$   \begin{array}{l} \det[sI_n-A]\\
                   \times \det[sI_q-F]\end{array}
    \times 
    \det\left[\begin{array}{cccc}
        I_n & 0 &        0        &-(sI_n-A)^{-1}B \\
         0  &I_q& (sI_q-F)^{-1}G  &        0       \\
         0  & 0 &       I_p       &C(sI_n-A)^{-1}B \\
         0  & 0 &H(sI_q-F)^{-1}G+K&      I_m    \end{array}\right] 
$$
And thus we obtain
\begin{eqnarray*}
  \varphi(s) &=&
    \det\left[\begin{array}{cc}
                I_p       &C(sI_n-A)^{-1}B \\
         H(sI_q-F)^{-1}G+K&      I_m    \end{array}\right]  \\
   && \quad\times \det[sI_n-A] \times \det[sI_q-F]\ .
\end{eqnarray*}

The off-diagonal entries in the first matrix are the transfer functions of the 
original system~(\ref{spls}) and of the compensator~(\ref{compensator}).
Consider  coprime factorizations
 \begin{eqnarray*}
  N(s)D(s)^{-1}&=&C(sI_n-A)^{-1}B\\
  P(s)Q(s)^{-1}&=&H(sI_q-F)^{-1}G+K\,.
 \end{eqnarray*}
Because $n$ and $q$ are the respective McMillan degrees, we have
 \begin{eqnarray*}
  \det D(s)&=& \det[sI_n-A]\,,\ \mbox{ and}\\
  \det Q(s)&=&\det[sI_q-F]\,,
 \end{eqnarray*}
and so our characteristic polynomial becomes
 \begin{equation}\label{closed-loop}
  \varphi(s) = \det\left[\begin{array}{cc}
                    Q(s) & N(s) \\
                    P(s) & D(s)  \end{array}\right]\ .
 \end{equation}
The first column of this 2 by 2 block matrix represents the Hermann-Martin
curve $M\colon{\mathbb P^1}\to\mbox{\it Grass}(p,{\mathbb K}^{m+p})$ of the
compensator and the second column the Hermann-Martin curve 
$L\colon{\mathbb P^1}\to\mbox{\it Grass}(m,{\mathbb K}^{m+p})$ of the
original system. 

The determinant (\ref{closed-loop}) must vanish at each root of the characteristic
polynomial.
Since, for every $s$, the columns giving the Hermann-Martin curves have full rank, 
we obtain the following version of the pole placement problem, when the
characteristic polynomial has distinct roots.
\medskip

\noindent{\bf Geometric Version of the pole placement problem.}
Suppose we have a strictly proper $m$-input $p$-output linear
system~(\ref{spls}) of McMillan degree $n$ with Hermann-Martin curve $L$
and a monic polynomial
$\varphi(s)$ of degree $n{+}q$ with distinct roots $s_1,s_2,\ldots,s_{n+q}$.
Which rational curves
$M\colon{\mathbb P^1}\to\mbox{\it Grass}(p,{\mathbb K}^{m+p})$
of degree $q$ satisfy
$$
  M(s_i)\cap L(s_i)\ \neq\ \{0\}\qquad\mbox{for}\quad  i=1,2,\ldots,n+q\ ?
$$
Thus we are looking for rational curves $M$ which satisfy $n+q$ Schubert
conditions of the type in Question 1.

Note that when $q=0$ (the case of static compensators),
$Q(s)=I_p$ and $P(s)=K$, so a static compensator is represented by the matrix
$$
   \left[\begin{array}{c}I_p\\K \end{array}\right]\,,
$$
whose column space is just a point in the Grassmannian
$\mbox{\it Grass}(p,{\mathbb K}^{m+p})$.
This observation of Byrnes~\cite{By80} was the point of departure for the
subsequent application of Schubert calculus to the pole placement problem.

\subsection{Number of dynamic compensators in the critical dimension}

Let ${\mathcal M}^q_{m,p}$ be the space of degree $q$ maps 
$M\colon{\mathbb P}^1\to\mbox{\it Grass}(p,{\mathbb K}^{m+p})$,
which is also the space of Hermann-Martin curves of possible
degree $q$ dynamic compensators for $m$-input, $p$-output linear
systems~(\ref{spls}). 
(This includes both proper and improper compensators.)
An important geometric perspective on the characteristic
equation~(\ref{closed-loop}) is that a given strictly proper linear system
of McMillan degree $n$ (represented by its Hermann-Martin curve
$L\colon{\mathbb P}^1\to\mbox{\it Grass}(p,{\mathbb K}^{m+p})$)
determines a {\it pole placement map}
$$
  \Lambda_L\ \colon\ {\mathcal M}^q_{m,p}\longrightarrow
\ \{\mbox{Polynomials of degree $n+q$}\}
$$
by
$$
  \Lambda_L\ \colon\ M\ \longmapsto\ \det[M(s) : L(s)]\,.
$$

This map gets its name from the fact that $\Lambda_L^{-1}(\varphi(s))$
is the set of Hermann-Martin curves of degree $q$ dynamic compensators
giving characteristic polynomial $\varphi(s)$.
Thus a strictly proper linear system is arbitrarily pole assignable
when the corresponding pole placement map is surjective.

Consider expanding this determinant along the columns of $M(s)$:
$$
  \Lambda_L(M)\ =\ \sum_{\alpha\in\binom{[m+p]}{p}}
                M_\alpha(s)\cdot L_\alpha(s)\,.
$$
Here $\binom{[m+p]}{p}$ is the collection of subsets of 
$\{1,2,\ldots,m{+}p\}$ of size $p$, $M_\alpha(s)$ is the
$\alpha$th maximal minor of $M(s)$ (given by the rows of $M(s)$ indexed by
$\alpha$), and $L_\alpha(s)$ is the appropriately 
signed complementary maximal minor of $L(s)$.
The point of this exercise is that the pole placement map is a linear
function of the coefficients of the polynomials $M_\alpha(s)$.

Thus we are led to consider the Pl\"ucker map
\begin{equation}\label{plucker}
  {\mathcal M}^q_{m,p}\ \longrightarrow\ 
   {\mathbb P}(\wedge^p{\mathbb K}^{m+p}\otimes{\mathbb K}^{q+1})
\end{equation}
which associates a $m{+}p$ by $p$ matrix $M(s)$ of polynomials 
(representing a degree $q$ compensator or degree $q$ curve) 
to its $\binom{m+p}{p}$ maximal minors $M_\alpha(s)$, which are
polynomials of degree $q$.
A more intrinsic definition of this map is given in Section~\ref{sec:reality}
just before~(\ref{eq:pluckermap}).
This gives a map to projective space as multiplying $M(s)$ by an invertible
$p$ by $p$ matrix $F$ multiplies each minor by the factor $\det F$ but
does not change the curve.
This Pl\"ucker map is an embedding, and one compactification of 
${\mathcal M}^q_{m,p}$ is the closure ${\mathcal K}^q_{m,p}$ of the image,
which we call the  {\it quantum Grassmannian}.
This space was introduced to systems theory by Rosenthal~\cite{Rosen_PhD}.

In this way, we see that the pole placement map factors
$$
  {\mathcal M}^q_{m,p}\ \longrightarrow\ 
  {\mathcal K}^q_{m,p}\ \stackrel{\pi_L}{\relbar\joinrel\longrightarrow}\ 
  {\mathbb P}^{n+q}\,,
$$
with the last map $\pi_L$ a linear projection on
${\mathbb P}(\bigwedge^p{\mathbb K}^{m+p}\otimes{\mathbb K}^{q+1})$.
Here, ${\mathbb P}^{n+q}$ is the space of polynomials of degree
at most $n+q$, modulo scalars.
(If the compensator is on the boundary of the compactification, then the
polynomial has degree less than $n{+}q$.)

Thus a necessary condition for arbitrary pole assignability of a 
strictly proper linear system $L$ is that 
$\pi_L$ is surjective.
The surjectivity of $\pi_L$  is sufficient for solving the pole placement
problem for $L$ and for generic polynomials $\varphi(s)$.
Rosenthal~\cite{Ro94} shows that if $q(m+p)+mp\leq n+q$ and
${\mathbb K}$ is algebraically closed, then $\pi_L$ is
surjective for generic strictly proper linear systems $L$.
This gives the criterion $n+q\leq q(m+p) + mp$ for a generic $m$-input, $p$-output
system~(\ref{spls}) of degree $n$ to be arbitrarily pole assigned
with degree $q$ compensators.

For generic systems $L$ in the critical dimension 
($q(m+p)+mp=n+q$ so that $\dim{\mathcal K}^q_{m,p}=n+q$)
the map $\pi_L$ is finite and hence surjective,
again, when ${\mathbb K}$ is algebraically closed.
Thus $\#(\pi_L^{-1}\varphi(s)\cap{\mathcal M}^q_{m,p})$ is the number of
compensators solving the pole placement problem for $\varphi(s)$.
Since $\pi_L$ is a linear projection, this number is bounded by 
the degree of the quantum Grassmannian ${\mathcal K}^q_{m,p}$ in its
Pl\"ucker embedding.
Suppose ${\mathbb K}$ is algebraically closed.
Since ${\mathcal M}^q_{m,p}$ is open in ${\mathcal K}^q_{m,p}$, for generic
$\varphi(s)$ the degree of ${\mathcal K}^q_{m,p}$ equals the number of
dynamic compensators, possibly counted with multiplicity.
This is the main theorem of~\cite{Ro94}.

When ${\mathbb K}={\mathbb R}$ so that $A,B,C$, and $\varphi(s)$ are real, 
$\pi^{-1}\varphi(s)\cap{\mathcal M}^q_{m,p}$ gives the 
{\it complex} dynamic compensators which solve the pole placement problem
for these data.
If $n+q\leq q(m+p)+mp$ and ${\mathcal K}^q_{m,p}$ has odd degree, then the set
of dynamic compensators is a projective variety defined over the real numbers of
odd degree, and hence contains a real point.
We deduce the following result.

\begin{theorem}\label{thm:odd}
Suppose $n\leq q(m+p-1)+mp$, and $\deg{\mathcal K}^q_{m,p}$ is odd.
Then a general strictly proper real linear system~(\ref{spls}) with $m$ inputs,
$p$ outputs, and McMillan degree $n$ is arbitrarily
pole assignable by real degree $q$ dynamic compensators.
\end{theorem}

When the degree of ${\mathcal K}^q_{m,p}$ is even the strongest result is
due to Rosenthal and Wang.

\begin{theorem}[\cite{RW96}]\label{thm:RW}
A generic  strictly proper linear system~(\ref{spls}) with $m$ inputs, $p$
outputs, and McMillan degree $n$ is arbitrarily pole assignable by
real degree $q$ compensators if 
$$
  n\ <\ q(m+p-1)+mp -\min\{ r_m(p-1), r_p(m-1)\}\,,
$$
where $r_p$ and $r_m$ are the remainders of $q$ upon division by $p$ and
$m$, respectively. 
\end{theorem}

The special case when $q=0$ of static compensation has an interesting history
(see the excellent survey of Byrnes~\cite{By89}).
In this case, the Grassmannian $\mbox{\it Grass}(p,{\mathbb K}^{m+p})$
plays the r\^{o}le of ${\mathcal K}^q_{m,p}$ and once it was discovered that
the equations for pole placement were linear equations on the Grassmannian in
its Pl\"ucker embedding, significant progress was made.
This included Brockett and Byrnes' calculation of the number of static
compensators for a generic $m$-input $p$-output linear system of McMillan
degree $mp$ as the degree of the Grassmannian~\cite{BB81}:
 \begin{equation}\label{eq:g-deg}
  (mp)!\frac{\prod_{1\leq j<k\leq p} (k-j)}{\prod_{j=1}^p(m+j-1)!}\ =\ 
 (mp)!\frac{1! \, 2! \, 3! \cdots (p\!- \!2) ! \, (p \!-\!1)!}{m!\, 
  (m \! + \! 1)!  \cdots(m \! + \! p \! - \! 1)!} 
 \end{equation}
We can deduce the analog of Theorem~\ref{thm:odd} from this;
unfortunately, this number is odd only when $\min(m,p)=1$ (and then it is 1),
or else $\min(m,p)=2$ and $\max(m,p)+1$ is a power of 2~\cite{Be76}.
The analog of Theorem~\ref{thm:RW} is due to Wang~\cite{Wa92}: $n<mp$ is
sufficient to guarantee arbitrary pole assignability over ${\mathbb R}$,
for generic systems.

\subsection{Formulae for $\deg{\mathcal K}^q_{m,p}$}\label{sec:Formulae}
\begin{figure}[thb]
$$\epsfxsize=5.52in \epsfbox{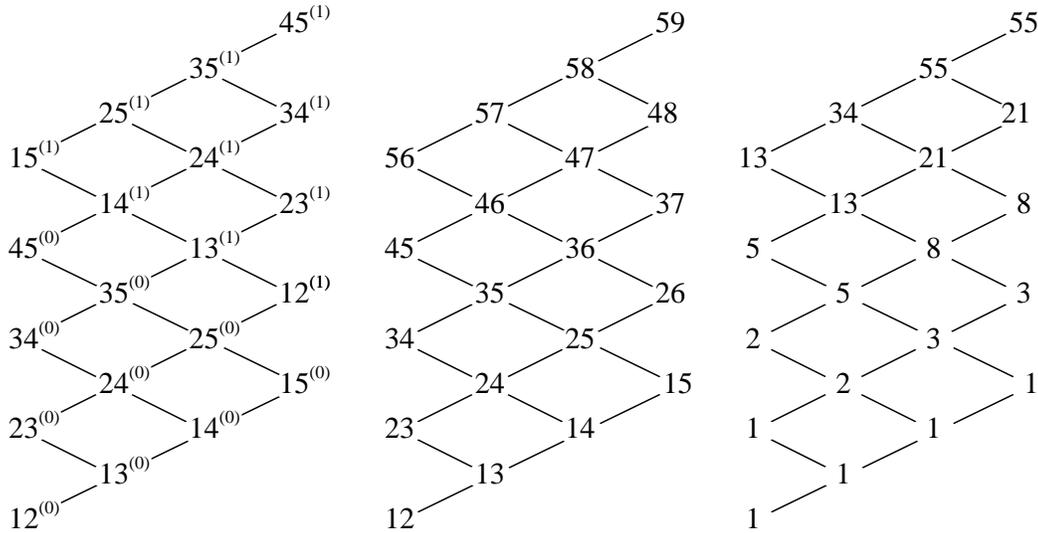}$$
 \caption{${\mathcal C}^1_{3,2}$, $J({\mathcal C}^1_{3,2})$, and 
           $\deg Z_{\alpha^{(a)}}$\label{fig:one}.}
\end{figure}

Let $z_{\alpha^{(a)}}$ be the coefficient of $s^a$ in the $\alpha$th maximal
minor $M_\alpha(s)$ of $M(s)$.
These coefficients provide {\it quantum Pl\"ucker coordinates} for 
${\mathbb P}(\bigwedge^p{\mathbb K}^{m+p}\otimes{\mathbb K}^{q+1})$.
Let 
${\mathcal C}^q_{m,p}:=
\{\alpha^{(a)}\mid \alpha\in\binom{[m+p]}{p}\mbox{ and }0\leq a\leq q\}$
be the indices of these quantum Pl\"ucker coordinates.
This index set has a natural partial order 
$$
  \alpha^{(a)}\ \leq\ \beta^{(b)} \quad\Longleftrightarrow\quad
   a\leq b \,\, \mbox{ and }  \,\,
   \alpha_i\leq \beta_{b-a+i} \,
  \mbox{ for }\,\, i = 1,2,\ldots,p-b+a\,. 
$$
The poset ${\mathcal C}^q_{m,p}$ is graded with the rank, $|\alpha^{(a)}|$,
of $\alpha^{(a)}$ equal to $a(m{+}p) + \sum_i(\alpha_i-i)$.
It is also a distributive lattice.
Figure~\ref{fig:one} shows ${\mathcal C}^1_{3,2}$ on the left.
Given $\alpha^{(a)}\in {\mathcal C}^q_{m,p}$, define the {\it quantum
Schubert variety}
\begin{equation}\label{qsv}
  Z_{\alpha^{(a)}}\quad :=\quad \{z\in {\mathcal K}^q_{m,p}\mid
  z_{\beta^{(b)}}=0\ \mbox{ if }\ \beta^{(b)}\not\leq \alpha^{(a)} \}\ .
\end{equation}
From this definition, we see that
$Z_{\alpha^{(a)}}\cap Z_{\beta^{(b)}}=Z_{\gamma^{(c)}}$ (set-theoretically),
where $\gamma^{(c)}$ is the greatest lower bound of $\alpha^{(a)}$ and
$\beta^{(b)}$.

Let ${\mathcal H}_{\alpha^{(a)}}$ be the hyperplane 
defined by $z_{\alpha^{(a)}}=0$.
We write $\beta^{(b)}\lessdot\alpha^{(a)}$ to indicate that
$\beta^{(b)}<\alpha^{(a)}$ and that there is no other
index $\gamma^{(c)}$ with 
$\beta^{(b)}<\gamma^{(c)}<\alpha^{(a)}$.
The main technical lemma of~\cite{RRW96,RRW98} is the following 

\begin{prop}[\cite{RRW96,RRW98}]\label{prop:RRW}
Let $\alpha^{(a)}\in{\mathcal C}^q_{m,p}$.
Then
\begin{enumerate}
\item[(i)]
   $Z_{\alpha^{(a)}}$ is an irreducible subvariety of 
   $\,{\mathcal K}^q_{m,p}$ of dimension $|\alpha^{(a)}|$.
\item[(ii)]
   The intersection of $\,Z_{\alpha^{(a)}}$ and 
   ${\mathcal H}_{\alpha^{(a)}}$ is
   generically transverse and we have 
   $$
     Z_{\alpha^{(a)}}\cap {\mathcal H}_{\alpha^{(a)}}\ =\ 
     \bigcup_{\beta^{(b)}\lessdot\alpha^{(a)}} Z_{\beta^{(b)}}\,.
   $$
\end{enumerate}
\end{prop}

This result is proven essentially by working in local coordinates
for $Z_{\alpha^{(a)}}$. 
Part (ii) is the geometric version of the (codimension-1) Pieri formula.
It generalizes the result of Schubert~\cite{Sch1886b}, 
who proved it for the classical Grassmannian ($a=0$).
By B\'ezout's Theorem (see~\cite[\S8]{Fu84a}), we deduce the following
fundamental recursion
 \begin{equation}\label{recursion}
    \deg Z_{\alpha^{(a)}}\ =\ 
     \sum_{\beta^{(b)}\lessdot\alpha^{(a)}} \deg Z_{\beta^{(b)}}\,.
 \end{equation}
The minimal quantum Schubert variety is a point, so we deduce a 
formula for $\deg{\mathcal K}^q_{m,p}$.

\begin{theorem}[\cite{RRW96,RRW98}]\label{degree}
  The degree $d(m,p;q)$ of $\,{\mathcal K}^q_{m,p}$  is the number of
  maximal chains in the poset ${\mathcal C}^q_{m,p}$.
\end{theorem}

For example, the degree of ${\mathcal K}^1_{3,2}$ is $55$, as shown by the
diagram on the right in Figure~\ref{fig:one}, which recursively computes the
degrees of the quantum Schubert varieties $Z_{\alpha^{(a)}}$.
In Section~\ref{sec:eqs} we give an alternative proof of Theorem~\ref{degree}
using Gr\"obner bases.

Ravi, Rosenthal, and Wang also solve this recursion to obtain the closed
formula~(\ref{formula}).
A first step is to change the indexing of the quantum Pl\"ucker coordinates,
embedding ${\mathcal C}^q_{m,p}$ into the set of increasing sequences 
$0<i_1<i_2<\cdots<i_p$ of
positive integers of length $p$.
Given   $\alpha^{(a)} \in {\mathcal C}^q_{p,m}$, write
$a=pl+r$ with $p > r \geq 0$ and define
a sequence $J(\alpha^{(a)})$ by
\begin{equation}\label{eq:seq_def}
  J(\alpha^{(a)})_k \ :=\ \left\{\begin{array}{ll}
   l(m{+}p)+\alpha_{r+k}&\quad\mbox{if }\ 1 \leq k\leq p-r\\
   (l{+}1)(m{+}p)+\alpha_{k-p+r}&\quad\mbox{if }\ p-r<k\leq p
  \end{array}\right.\,.
\end{equation}
For instance, when $m=p=5$, we have
$J((2,3,5,6,9)^{(7)})=(15,16,19,22,23)$.
Note that we have $J(\alpha^{(a)})_p<J(\alpha^{(a)})_1+m+p$.
This gives an order isomorphism of 
the poset ${\mathcal C}^q_{p,m}$ with  
the poset of sequences $i_1<i_2<\cdots<i_p$ of positive integers where 
$i_p<i_1+m+p$.
This is illustrated in the middle diagram of Figure~\ref{fig:one}, which shows
the image of ${\mathcal C}^1_{3,2}$.
This isomorphism (of course) preserves the rank function of the two 
posets:
 \begin{equation}\label{Miracle}
  | \alpha^{(a)} | \ := \ 
  a(m+p) + \sum_{j=1}^p (\alpha_j - j) 
  \ = \ 
  \sum_{i=1}^p \bigl(  J(\alpha^{(a)})_i - i \,\bigr)
  \ =: \ 
  | J(\alpha^{(a)})|\ .
 \end{equation}
Observe that $J(\alpha^{(a)})$ is congruent to $\alpha$ modulo $m+p$.

\begin{lemma}\label{lem:recursion}
Let $d(i_1,i_2,\ldots,i_p)$ be a function defined for all 
weakly increasing sequences of non-negative integers $i_1,i_2,\ldots,i_p$
with $i_p\leq i_1+m+p$.
Suppose that for any sequence $0<i_1<\cdots<i_p$ with $i_p<i_1{+}m{+}p$, 
$d(i_1,i_2,\ldots,i_p)$ satisfies the recursion 
 \begin{equation}\label{eq:rec}
  d(i_1,i_2,\ldots,i_p)\ =\ \sum_{k=1}^p 
   d(i_1,i_2,\ldots,i_k-1,\ldots,i_p)\,,
 \end{equation}
is subject to the initial condition 
 \begin{equation}\label{initial}
    d(1,2,\ldots,p)=1\,,
 \end{equation}
and the boundary conditions 
 \begin{eqnarray}
  d(\ldots,l,l,\ldots)&=&0\,\makebox[1em][l]{,} \label{A}\\
  d(0,\ldots) &=& 0\,\makebox[1em][l]{,}\mbox{ and}\label{B}\\
  d(i_1,i_2,\ldots,i_p) &=& 0\,\makebox[1em][l]{\ }\mbox{ if }\ 
        i_p=i_1+m+p\label{C}\,.
 \end{eqnarray}
Then $d(J(\alpha^{(a)}))=\deg Z_{\alpha^{(a)}}$.
\end{lemma}

\begin{proof}
Let $j_1,j_2,\ldots,j_p=J(\alpha^{(a)})$ for 
$\alpha^{(a)}\in{\mathcal C}^q_{m,p}$.
Then the sequence $I:=j_1,j_2,\ldots,j_k{-}1$, $\ldots,j_p$ fails to equal 
$J(\beta^{(b)})$ for some $\beta^{(b)}\in{\mathcal C}^q_{m,p}$
only if $I$ has either two repeated indices~(\ref{A}), or if
$i_1=0$~(\ref{B}), or else if $i_p=i_1+m+p$~(\ref{C}).
In each of these cases $d(I)=0$, and so the function $d(J(\alpha^{(a)}))$
defined for $\alpha^{(a)}\in{\mathcal C}^q_{m,p}$ satisfies the
recursion~(\ref{recursion}) for $\deg Z_{\alpha^{(a)}}$.

Since the index of a minimal quantum Schubert variety (which is a point)
is $(1,2,\ldots,p)^{(0)}$, and $J(\alpha^{(0)})=\alpha$, the function
$d(J(\alpha^{(a)}))$ also satisfies the initial condition for 
$\deg Z_{\alpha^{(a)}}$.
\end{proof}

Sequences $I:0<i_1<i_2<\cdots<i_p\leq l$ index Schubert varieties
$\Omega_I$ of $\mbox{\it Grass}(p,{\mathbb K}^l)$.
Schubert~\cite{Sch1886b} showed that  the degree $g(I)$ 
of the Schubert variety $\Omega_I$ in the
Pl\"ucker embedding satisfies the recursion, initial condition, and
boundary conditions~(\ref{A}) and~(\ref{B}) of Lemma~\ref{lem:recursion}.
He later~\cite{Sch1886c} gave the following closed formula for this degree
(compare with~(\ref{eq:g-deg})):
 \begin{equation}\label{Sch_formula}
  |I|!\frac{\prod_{j<k}(i_k-i_j)}{\prod_j(i_j-1)!}\,,
 \end{equation}
where $|I|=\sum_j i_j-j$.
This formula~(\ref{Sch_formula}) defines $g(I)$ as an alternating function on
all sequences of integers if we set $1/l!=0$ when $l<0$.

\begin{theorem}[\cite{RRW96}]\label{sum}
  Let $\alpha^{(a)}\in{\mathcal C}^q_{m,p}$ and set $I:=J(\alpha^{(a)})$.
  Then we have
  $$
    d(I)\ =\ \sum_{b_1+\cdots+b_p=0}
             g(i_1+b_1(m+p),\, i_2+b_2(m+p),\,\ldots,\,
               i_p+b_p(m+p))\,.
  $$
\end{theorem}

Observe that the sum is in fact finite, as only sequences
$b_1,b_2,\ldots,b_m$ for which every term 
$i_j+b_j(m+p)$ is positive contribute to the sum.

\begin{proof}
Let $\delta(I)$ be the function defined by the sum.
First observe that if $\max I \leq m+p$, then there is only the trivial
summand (all $b_i=0$) and so $\delta(I)=g(I)$.
Also, since $g$ is alternating, $\delta$ is an alternating function.

We show that the function $\delta$ satisfies the conditions of
Lemma~\ref{lem:recursion}, when $I$ is a weakly increasing sequence of
non-negative integers with $i_p\leq i_1+m+p$.
First, $\delta$ satisfies the recursion of Lemma~\ref{lem:recursion}
because the function $g$ satisfies the recursion.
Second, $\delta(1,2,\ldots,p)=g(1,2,\ldots,p)=1$, giving the initial
condition. 
Next, since $\delta$ is alternating, it
satisfies~(\ref{A}).
Suppose $i_p=i_1+m+p$.
Then every summand indexed by $b_1,b_2,\ldots,b_p$ with $b_p-1=b_1$ vanishes as
$g$ is alternating, and every summand with $b_p{-}1\neq b_1$ is paired with
another summand indexed by $b_p{+}1,b_2,\ldots,b_1{-}1$, which has the same
absolute value, but opposite sign, as $g$ is alternating.
Thus~(\ref{C}) holds for $\delta$.
Finally, if $i_1=0$, then either $i_p=m+p$ and so $\delta(I)=0$ or else
$i_p<0$ and so $\delta(I)=g(I)=0$, giving~(\ref{B})
and proving the theorem.
\end{proof}

We now deduce the formula~(\ref{formula}) from these results.
The quantum Grassmannian ${\mathcal K}^q_{m,p}$ is the maximal
quantum Schubert variety $Z_{(m+1,\ldots,m+p)^{(q)}}$.
Let $q=pl+r$ with $0\leq r<p$.
Set $\alpha:=(m{+}1,\ldots,m{+}p)$ and $n:=m+p$.
Then $J(\alpha^{(q)})$ is the sequence
$$
  (ln+m+r+1,\,\ldots,\,ln+m+p,\,(l+1)n+m+1,\,\ldots,
     (l+1)n+m+r)
$$
and we have
$$
  |\alpha^{(q)}|\ =\ |J(\alpha^{(q)})|\ =\ 
   pln+rn+mp\ =\ q(m{+}p)+mp\,.
$$

By Theorem~\ref{sum}, the degree $d(J(\alpha^{(q)}))$ of 
${\mathcal K}^q_{m,p}$ is the sum over all sequences of integers
$b_1,b_2,\ldots,b_p$ satisfying
$b_1+b_2+\cdots+b_p=0$
of the terms\medskip

\noindent\qquad$g((l+b_1)n+m+r+1,\,\ldots,\,(l+b_{p-r})n+m+p,$ \smallskip

\hfill$(l+b_{p-r+1}+1)n+m+1,\,\ldots,\,(l+b_p+1)n+m+r)$.\qquad\medskip

This term equals
$$
  (-1)^{r(p-r)}\,g(\nu_1n+m+1,\,\nu_2n+m+2,\,\ldots,\,\nu_pn+m+p)\,,
$$
where 
$$
  (\nu_1,\nu_2,\ldots,\nu_p)\ =\ 
  (l+b_{p-r+1}+1,\ldots,l+b_p+1,\,l+b_1,\ldots,l+b_{p-r})\,,
$$
and the sign $(-1)^{r(p-r)}$ comes from the resulting permutation of the
argument of $g$.
Since $\nu_1+\cdots+\nu_p=q$, we obtain the formula
$$
  \deg{\mathcal K}^q_{m,p}\ =\ 
    (-1)^{r(p-r)}\,\sum_{\nu_1+\cdots+\nu_p=q}
    g(\nu_1n+m+1,\nu_2n+m+2,\ldots,\nu_pn+m+p)\,.
$$
Finally, to obtain the formula~(\ref{formula}), we use Schubert's 
formula~(\ref{Sch_formula}) for $g$ and the following identity
\begin{eqnarray*}
  r(p-r)&=&(q-pl)((l+1)p-q)\ =\ 2pql+pq-q^2-p^2l(l+1)\\
        &\equiv& pq+q^2\mod 2\ \equiv\ pq+q\mod 2\,.
\end{eqnarray*}
%


\section{Reality and Transversality}\label{sec:reality}
Traditionally, intersection theory and enumerative geometry (both classical
and quantum) treat the case of complex solutions to enumerative problems,
for it is in this case that the most general and elegant results hold.
The real numbers pose special problems as the number of real figures
satisfying conditions imposed by general (fixed) real figures depends subtly
on the configuration of the fixed real figures.
Algebraically closed fields of positive characteristic also 
pose special problems
in enumerative geometry as the number of solutions may depend upon the
characteristic of the field.
One reason for this is that the solutions may occur with multiplicities; 
the subvarieties defined by the conditions may not intersect transversally
in positive characteristic.
In characteristic zero, Kleiman's Theorem on the 
transversality of a general translate~\cite{MR50:13063} may be invoked to
show that each solution to many enumerative problems (including that of
Question 1) occurs without multiplicities. 
In positive characteristic, general translates are not necessarily
transverse, and other techniques must be employed to determine whether 
the solutions occur without multiplicity.

For the enumerative problem of Question 1, both these difficulties may be
overcome using the same elementary arguments, which are a version of the
theory of~\cite{So_flags} adapted to this particular enumerative problem.
These arguments are based upon the Pieri homotopy algorithm of~\cite{HSS}
and related to a numerical homotopy continuation algorithm for computing
numerical solutions to these enumerative problems when 
${\mathbb K}={\mathbb C}$.
See~\cite[\S5]{HV} for details of this algorithm.

\begin{theorem}[\cite{So00,So_trans}]\label{real-trans}
Let $m,p>1$ and $q\geq 0$ be integers.
Set $N:=q(m+p)+mp$.
Suppose ${\mathbb K}$ is an infinite field with algebraic closure
$\overline{\mathbb K}$.
Then there exist points
$s_1,s_2,\ldots,s_N\in{\mathbb P}^1_{\mathbb K}$ and $m$-planes
$L_1,L_2,\ldots,L_N\subset{\mathbb K}^{m+p}$ for which there are
exactly $d(m,p;q)$ maps 
$M\colon{\mathbb P}^1\to\mbox{\it Grass}(p,\overline{\mathbb K}^{m+p})$ 
of degree $q$ satisfying
$$
  M(s_i)\cap L_i\ \neq\ \{0\}\quad\mbox{for } i=1,\ldots,N\,.
$$

When ${\mathbb K}={\mathbb R}$, we may further choose the real points and real $m$-planes 
so that all of the resulting maps are real.
\end{theorem}

Thus the enumerative problem of Question 1 is enumerative in all
characteristics
and when ${\mathbb K}={\mathbb R}$, there is some choice of points and $m$
planes for which all of the {\it a priori} complex solutions are
real.\smallskip

Suppose ${\mathbb K}$ is an infinite field.
Let $L\subset{\mathbb K}^{m+p}$ be an $m$-plane, none of whose Pl\"ucker
coordinates vanish.
That is, if $L$ is the column space of a $m+p$ by $m$ matrix, also written
$L$, then none of the $m$ by $m$ maximal minors of $L$ vanishes.
This choice is possible as ${\mathbb K}$ is infinite.
Let ${\bf e}_1,{\bf e}_2,\ldots,{\bf e}_{m+p}$ be the distinguished basis of
${\mathbb K}^{m+p}$ 
corresponding to the rows of such matrices.
We equip ${\mathbb K}^{m+p}$ with an action of 
${\mathbb K}^\times$.
For $s\in{\mathbb K}^\times$, set
$$
  (s,{\bf e}_i)\ \longmapsto\ s.{\bf e}_i\ :=\ s^{m+p-i}{\bf e}_i\,.
$$

For $s\in{\mathbb K}^\times$, 
let ${\mathcal H}(s,L)$ be the hyperplane in Pl\"ucker space
${\mathbb P}(\bigwedge^p{\mathbb K}^{m+p}\otimes{\mathbb K}^{q+1})$
whose intersection with the space ${\mathcal M}^q_{m,p}$ of rational curves
of degree $q$ consists of those maps 
$M\colon{\mathbb P}\to\mbox{\it Grass}(p,{\mathbb K}^{m+p})$ satisfying
$$
  M(s^{m+p})\cap s.L\ \neq\ \{0\}\,.
$$
This is just the set of maps satisfying 
$M(t)\cap K\neq\{0\}$, where $t=s^{m+p}$ and $K=s.L$.

This condition is equivalent to the vanishing of the determinant
$$
  \det [M(s^{m+p}): s.L]\,.
$$
If we expand this determinant along the columns of $M(s^{m+p})$, we obtain
$$
  \sum_\alpha M_\alpha(s^{m+p})\,L_\alpha\,s^{\binom{m}{2}}\,s^{|\alpha|}\,,
$$
where $M_\alpha(s)$ is the $\alpha$th maximal minor of $M(s)$ and 
$L_\alpha$ is the appropriately signed complementary maximal minor
(Pl\"ucker coordinate) of $L$.
If we now expand the polynomials $M_\alpha(s^{m+p})$ in terms of the quantum
Pl\"ucker coordinates $z_{\alpha^{(a)}}$ of $M$ and divide out the common
factor $s^{\binom{m}{2}}$, 
we obtain the following equation for  
the hyperplane ${\mathcal H}(s,L)$:
$$
  \Phi(s,L)\ =\ \sum_{\alpha^{(a)}\in{\mathcal C}^q_{m,p}}
                z_{\alpha^{(a)}}\,L_\alpha\,s^{|\alpha^{(a)}|}\,,
$$
since $a(m+p) + |\alpha| = |\alpha^{(a)}|$.

We prove Theorem~\ref{real-trans} by first showing that given any
$m$-plane $L\subset{\mathbb K}^{m+p}$ with no vanishing
Pl\"ucker coordinates, 
there exist points $s_1,s_2,\ldots,s_N\in{\mathbb K}$ 
such that there are exactly $d(m,p;q)$ points in 
${\mathcal K}^q_{m,p}$ common to the hyperplanes
${\mathcal H}(s_i,L)$, and then argue that none of these
$d(m,p;q)$ points lie in the boundary
${\mathcal K}^q_{m,p}-{\mathcal M}^q_{m,p}$.

We first study the boundary of ${\mathcal K}^q_{m,p}$, using results of
Bertram~\cite{Be97}. 
While Bertram works over the complex numbers ${\mathbb C}$, his results we
invoke remain valid over any field.
A smooth compactification of ${\mathcal M}^q_{m,p}$ is provided by a quot scheme 
${\mathcal Q}^q_{m,p}$~\cite{St87}.
By definition, there is a universal exact sequence
$$
  0\ \rightarrow\ {\mathcal S}\ \rightarrow\
  {\mathbb K}^{m+p}\otimes
  {\mathcal O}_{{\mathbb P}^1\times{\mathcal Q}^q_{m,p}}\ 
  \rightarrow\  {\mathcal T}\ \rightarrow 0
$$
of sheaves on ${\mathbb P}^1\times{\mathcal Q}^q_{m,p}$ where 
${\mathcal S}$ is a vector bundle of degree $-q$ and rank $p$.
Twisting the determinant of ${\mathcal S}$ by 
${\mathcal O}_{{\mathbb P}^1}(q)$ and pushing forward to 
${\mathcal Q}^q_{m,p}$ induces a Pl\"ucker map
 \begin{equation}\label{eq:pluckermap}
  {\mathcal Q}^q_{m,p}\ \rightarrow\ 
  {\textstyle {\mathbb P}\left( \bigwedge^p{\mathbb K}^{m+p}\otimes 
      H^0({\mathcal O}_{{\mathbb P}^1}(q))^*\right)}\,.
 \end{equation}
The restriction to ${\mathcal M}^q_{m,p}$ is the Pl\"ucker
map~(\ref{plucker}) and the image is the quantum Grassmannian
${\mathcal K}^q_{m,p}$.

The Pl\"ucker map fails to be injective on the boundary 
${\mathcal Q}^q_{m,p}-{\mathcal M}^q_{m,p}$ of ${\mathcal Q}^q_{m,p}$.
Indeed, Bertram constructs a ${\mathbb P}^{p-1}$ bundle over 
${\mathbb P}^1\times{\mathcal Q}^{q-1}_{m,p}$ that maps onto the boundary of  
${\mathcal Q}^q_{m,p}$, with its
restriction over ${\mathbb P}^1\times{\mathcal M}^{q-1}_{m,p}$ an embedding.
On this projective bundle, the Pl\"ucker map factors through the base 
${\mathbb P}^1\times{\mathcal Q}^{q-1}_{m,p}$ and the image of a point in
the base is $s\cdot M$, where $s$ is the section of 
${\mathcal O}_{{\mathbb P}^1}(1)$ vanishing at $s\in{\mathbb P}^1$ and 
$M$ is the image of a point in ${\mathcal Q}^{q-1}_{m,p}$ under its
Pl\"ucker map. 

This identifies the image of the exceptional locus of the Pl\"ucker map,
which is the boundary of ${\mathcal K}^q_{m,p}$,  
with the image of ${\mathbb P}^1\times{\mathcal K}^{q-1}_{m,p}$ in 
${\mathcal K}^q_{m,p}$ under a map $\pi$ which we now describe.
Let $x_{\alpha^{(a)}}$ be the quantum Pl\"ucker coordinates
for ${\mathcal K}^{q-1}_{m,p}$.
Then the boundary of ${\mathcal K}^q_{m,p}$ is the image of the map
$\pi:{\mathbb P}^1\times{\mathcal K}^{q-1}_{m,p}\rightarrow
         {\mathcal K}^q_{m,p}$
defined by 
\begin{equation}\label{eq:pi}
   \left([A,B],(x_{\beta^{(b)}} : 
   \beta^{(b)}\in{\mathcal C}^{q-1}_{m,p})\right)
   \ \longmapsto\ 
  (Ax_{\alpha^{(a)}}-Bx_{\alpha^{(a-1)}} :
      \alpha^{(a)}\in{\mathcal C}^q_{m,p})\,,
\end{equation}
where $x_{\alpha^{(q)}}=x_{\alpha^{(-1)}}=0$.

For a variety $X$ defined over ${\mathbb K}$, let 
$X(\overline{\mathbb K})$ be the $\overline{\mathbb K}$-points of $X$.
Theorem~\ref{real-trans} is a consequence of the following two theorems.

\begin{theorem}\label{thm:transverse}
Suppose $L\subset{\mathbb K}^{m+p}$ is an 
$m$-plane with no vanishing Pl\"ucker coordinates.
Then there exist $s_1,s_2,\ldots,s_N\in{\mathbb K}$ so that the intersection 
\begin{equation}\label{thm31}
  Z_{\alpha^{(a)}}(\overline{\mathbb K})\cap{\mathcal H}(s_1,L)
  \cap{\mathcal H}(s_2,L)\cap\cdots
  \cap{\mathcal H}(s_{|\alpha^{(a)}|},L)
\end{equation} 
is transverse for any $\alpha^{(a)}\in{\mathcal C}^q_{m,p}$.

If\/ ${\mathbb K}={\mathbb R}$, then we may further choose these numbers 
$s_1,s_2,\ldots,s_N$ so that for any
$\alpha^{(a)}\in{\mathcal C}^q_{m,p}$, all points in the
intersection~(\ref{thm31}) are real.
\end{theorem}

\begin{theorem}\label{thm:proper}
Suppose $L\subset{\mathbb K}^{m+p}$ is an 
$m$-plane with no vanishing Pl\"ucker coordinates.
If $s_1,s_2,\ldots,s_k\in {\mathbb K}$ are distinct, 
then for any $\alpha^{(a)}\in{\mathcal C}^q_{m,p}$
the intersection  
 \begin{equation}\label{eq:proper}
   Z_{\alpha^{(a)}}\cap{\mathcal H}(s_1,L)
   \cap{\mathcal H}(s_2,L)\cap\cdots\cap{\mathcal H}(s_k,L)
 \end{equation}
is proper in that it has dimension $|\alpha^{(a)}|-k$.
\end{theorem}

\begin{proof}[Proof of Theorem~\ref{real-trans}]
By Theorem~\ref{thm:transverse}, there exist 
$s_1,s_2,\ldots,s_N\in{\mathbb K}$ so that the intersection 
\begin{equation}\label{eq:int}
  {\mathcal K}^q_{m,p}(\overline{\mathbb K}) \cap 
     {\mathcal H}(s_1,L)\cap
   {\mathcal H}(s_2,L)\cdots\cap{\mathcal H}(s_N,L)
\end{equation}
is transverse and consists of exactly $d(m,p;q)$ points,
and when ${\mathbb K}={\mathbb R}$, these points of intersection are real.
Furthermore, we may choose these numbers $s_i$ so that
their $(m{+}p)$th powers are distinct.
To prove Theorem~\ref{real-trans}, we show that these points all lie in 
${\mathcal M}^q_{m,p}$.
Thus each point in~(\ref{eq:int}) represents a map 
$M:{\mathbb P}^1\to\mbox{\it Grass}(p,{\mathbb K}^{m+p})$
of degree $q$ satisfying
$M(s_i^{m+p})\cap s_i.L\neq\{0\}$ for $i=1,2,\ldots,N$.

Let $\pi:{\mathbb P}^1\times{\mathcal K}^{q-1}_{m,p}\to{\mathcal K}^q_{m,p}$ 
be the map~(\ref{eq:pi}) whose
image is the complement of ${\mathcal M}_q$ in ${\mathcal K}_q$.
Then 
\begin{eqnarray*} 
  \pi^*\Phi(s,L)&=&
   \sum_{\alpha^{(a)}\in{\mathcal C}^q_{m,p}}
    (Ax_{\alpha^{(a)}}-Bx_{\alpha^{(a-1)}})\,
       L_\alpha\,s^{|\alpha^{(a)}|}\\
  &=&
    (A-Bs^{m+p})\sum_{\beta^{(b)}\in{\mathcal C}^{q-1}_{m,p}}
    x_{\beta^{(b)}}\,L_\beta\,s^{|\beta^{(a)}|}\\
  &=&
    (A-Bs^{m+p})\ \Phi'(s,L)\,,
\end{eqnarray*}
where $\Phi'(s,L)$ is the linear form for ${\mathcal K}^{q-1}_{m,p}$ analogous
to $\Phi(s,L)$.
Let ${\mathcal H}'(s,L)$ be the hyperplane given by the linear form
$\Phi'(s,L)$.

Any point in the intersection~(\ref{eq:int}) but not in ${\mathcal M}^q_{m,p}$ is
the image of a point 
$([A,B],x)$ in ${\mathbb P}^1\times{\mathcal K}^{q-1}_{m,p}$
satisfying 
$\pi^*\Phi(s_i,L)=(A-Bs_i^{m+p})\Phi'(s_i,L)$ 
for each $i=1,2,\ldots,N$.
As the $(m{+}p)$th powers of the $s_i$ are distinct, 
such a point can satisfy $A-Bs_i^{m+p}=0$ for at most one $i$.
Thus $x\in{\mathcal K}^{q-1}_{m,p}$ lies in 
at least $N{-}1$ of the hyperplanes ${\mathcal H}'(s_i,L)$.
Since $N{-}1$ exceeds the dimension $N{-}(m{+}p)$ of 
${\mathcal K}^{q-1}_{m,p}$, there are no such 
points $x\in{\mathcal K}^{q-1}_{m,p}$, by Theorem~\ref{thm:proper} 
applied to maps of degree $q{-}1$.
\end{proof}

\begin{proof}[Proof of Theorem~\ref{thm:proper}]
For any $s_1,s_2,\ldots,s_k$, the intersection~(\ref{eq:proper}) has dimension 
at least $|\alpha^{(a)}|-k$.
We show it has at most this dimension, if $s_1,s_2,\ldots,s_k$ are distinct.

Suppose $k=|\alpha^{(a)}|+1$ and let $z\in  Z_{\alpha^{(a)}}$.
Then $z_{\beta^{(b)}}=0$ if $\beta^{(b)}\not\leq\alpha^{(a)}$
and so the form $\Phi(s,L)$ defining ${\mathcal H}(s,L)$ 
evaluated at $z$ is 
$$
\sum_{\beta^{(b)}\leq\alpha^{(a)}} z_{\beta^{(b)}}\,L_\beta\,
s^{|\beta^{(b)}|}\,.
$$
This is a non-zero polynomial in $s$ of degree at most $|\alpha^{(a)}|$
and thus it vanishes for at most $|\alpha^{(a)}|$ distinct values of $s$.
It follows that~(\ref{eq:proper}) is empty for $k>|\alpha^{(a)}|$.

If $k\leq |\alpha^{(a)}|$ and $s_1,s_2,\ldots,s_k$ are distinct,
but~(\ref{eq:proper}) has dimension exceeding  
$|\alpha^{(a)}|-k$, then we may complete $s_1,s_2,\ldots,s_k$ to a set of distinct 
numbers $s_1,s_2,\ldots,s_{|\alpha^{(a)}|+1}$ which give a non-empty
intersection in~(\ref{eq:proper}), a contradiction.\end{proof}

\begin{proof}[Proof of Theorem~\ref{thm:transverse}]
We prove both parts of the theorem simultaneously, making note of
the differences when ${\mathbb K}={\mathbb R}$.

We construct the sequence $s_i$ inductively.
The unique element of rank 1 in ${\mathcal C}^q_{m,p}$ is $\alpha^{(0)}$,
where  $\alpha$ is the sequence $1<2<\cdots<p{-}1<p{+}1$.
The quantum Schubert variety $Z_{\alpha^{(0)}}$ is a line in Pl\"ucker space. 
Indeed, it is isomorphic to the set of $p$-planes containing a fixed
$(p{-}1)$-plane and lying in a fixed $(p{+}1)$-plane.
By Theorem~\ref{thm:proper} or direct observation, 
$Z_{\alpha^{(0)}}\cap{\mathcal H}(s,L)$ is
then a single point, for any non-zero $s$.
When ${\mathbb K}={\mathbb R}$, this point is real.
Let $s_1\in{\mathbb K}^\times$ be arbitrary.

Suppose $s_1,s_2,\ldots,s_k\in{\mathbb K}$ are distinct points with the property 
that for any $\beta^{(b)}$ with $|\beta^{(b)}|= k$, 
$$
  Z_{\beta^{(b)}}\cap{\mathcal H}(s_1,L)\cap{\mathcal H}(s_2,L)
   \cap\cdots\cap{\mathcal H}(s_{|\beta^{(b)}|},L)
$$ 
is transverse. 
When ${\mathbb K}={\mathbb R}$, we suppose further that all points of
intersection are real.

Let $\alpha^{(a)}$ be an index with $|\alpha^{(a)}|=k+1$ and consider the
1-parameter family ${\mathcal Z}(s)$ of schemes defined by 
$ Z_{\alpha^{(a)}}\cap{\mathcal H}(s,L)$, for $s\in{\mathbb K}^\times$.
If we restrict the form $\Phi(s,L)$ to 
$z\in Z_{\alpha^{(a)}}$, then
we obtain
$$
  \sum_{\beta^{(b)}\leq\alpha^{(a)}} 
  z_{\beta^{(b)}}\,L_\beta\, s^{|\beta^{(b)}|}\,,
$$
a polynomial in $s$ with leading term 
$z_{\alpha^{(a)}}\,L_\alpha\, s^{|\alpha^{(a)}|}$.
Since the Pl\"ucker coordinate $L_\alpha$ is non-zero, 
${\mathcal Z}(\infty)\subset Z_{\alpha^{(a)}}$ is defined by
$z_{\alpha^{(a)}}=0$, and so ${\mathcal Z}(\infty)$ equals
$$
    Z_{\alpha^{(a)}}\cap {\mathcal H}_{\alpha^{(a)}}\ =\ 
   \bigcup_{\beta^{(b)}\lessdot\alpha^{(a)}}  Z_{\beta^{(b)}}\,,
$$
by Proposition~\ref{prop:RRW} (ii).\smallskip

\noindent{\bf Claim:}
The cycle
$$
   {\mathcal Z}(\infty)
   \cap{\mathcal H}(s_1,L)\cap{\mathcal H}(s_2,L)
   \cap\cdots\cap{\mathcal H}(s_k,L)
$$
is free of multiplicities.

If not, then there are two components $ Z_{\beta^{(b)}}$
and $ Z_{\gamma^{(c)}}$ of ${\mathcal Z}(\infty)$ such that 
$$
  Z_{\beta^{(b)}} \cap Z_{\gamma^{(c)}}
  \cap {\mathcal H}(s_1,L)\cap{\mathcal H}(s_2,L)\cap 
  \cdots\cap {\mathcal H}(s_k,L)
$$
is non-empty.
But this contradicts Theorem~\ref{thm:proper}, as 
$Z_{\beta^{(b)}} \cap Z_{\gamma^{(c)}}=Z_{\delta^{(d)}}$, where
$\delta^{(d)}$ is the greatest lower bound of $\beta^{(b)}$
and $\gamma^{(c)}$ in ${\mathcal C}^q_{m,p}$, and so  
$\dim Z_{\delta^{(d)}}<\dim Z_{\beta^{(b)}}=k$.

Because the intersection of ${\mathcal Z}(\infty)$, the fibre of 
${\mathcal Z}$ at infinity, with the cycle 
${\mathcal H}(s_1,L)\cap{\mathcal H}(s_2,L)
        \cap\cdots\cap{\mathcal H}(s_k,L)$
is zero dimensional and free of multiplicities, it is transverse,
and so the general fibre of ${\mathcal Z}$ meets
${\mathcal H}(s_1,L)\cap{\mathcal H}(s_2,L)\cap\cdots\cap{\mathcal H}(s_k,L)$
transversally.
Thus there is a non-empty Zariski open subset $O_{\alpha^{(a)}}$ of 
${\mathbb A}^1_{\mathbb K}$ consisting of points $s$ for which 
$$
  Z_{\alpha^{(a)}}\cap{\mathcal H}(s_1,L)\cap{\mathcal H}(s_2,L)\cap
  \cdots\cap{\mathcal H}(s_k,L)\cap{\mathcal H}(s,L)
$$ 
is transverse.
Choose $s_{k+1}$ to be any point common to all $O_{\alpha^{(a)}}$ for
$|\alpha^{(a)}|=k+1$.

When ${\mathbb K}={\mathbb R}$, 
the claim implies there is a real number ${\mathcal N}_{\alpha^{(a)}}>0$ 
such that if $s>{\mathcal N}_{\alpha^{(a)}}$, then 
$$
  {\mathcal Z}(s)\cap{\mathcal H}(s_1,L)\cap\cdots\cap{\mathcal H}(s_k,L)
$$
is transverse with all points of intersection real.
Set
$$
  {\mathcal N}_{k+1}\ :=\ 
  \max\{{\mathcal N}_{\alpha^{(a)}} : |\alpha^{(a)}|=k+1\}\,
$$
and let $s_{k+1}$ be any real number satisfying
$s_{k+1}>{\mathcal N}_{k+1}$. 
\end{proof}

\begin{remark}
While these results rely upon work from systems theory, the result when
${\mathbb K}={\mathbb R}$ unfortunately does not give any insight into the dynamic pole
placement problem:
In the dynamic pole placement problem, the planes $L(s)$ lie on a rational curve
$L(s)$ of degree $mp+q(m+p)-q$ while the planes $s_i.L$ of
Theorem~\ref{thm:transverse} lie on the rational curve $s.L$, which has degree
$mp$. 
Thus there is overlap only when $q=0$, which is the static pole placement
problem.
\end{remark}

\begin{remark}
Theorem~\ref{thm:transverse} proves reality and transversality for the
enumerative problem of Question 1.
There are more general enumerative problems involving rational curves on a
Grassmannian obtained by replacing the Schubert condition
$M(s)\cap L\neq\{0\}$ with more general Schubert conditions.
It is not known, but is expected, that the transversality and reality properties
established in Theorem~\ref{thm:transverse} for the enumerative problem of
Question 1 hold also for these more general enumerative problems.
\end{remark}

\begin{remark}\label{rem:intersection}
The argument given at the end of the Proof of Theorem~\ref{real-trans},
that the intersection~(\ref{eq:int}) contains no points in the boundary and
hence lies in ${\mathcal M}^q_{m,p}$, may be generalized to show that when 
$a=q$, the intersection~(\ref{thm31}) similarly lies in 
${\mathcal M}^q_{m,p}$.
\end{remark}


\section{Equations for the Quantum Grassmannian}\label{sec:eqs}
In Section~\ref{sec:reality}, we solved the enumerative problem of Question 1 by arguing
directly from the equations describing the conditions~(\ref{condition}).
This is an unusual feature of that enumerative problem: despite the fact
that algebraic geometry is ostensibly concerned with solutions to
polynomial equations, enumerative geometric problems are not typically
solved in this manner.
What is more unusual is that this enumerative problems admits a second
solution also based upon equations, in this case equations for the quantum
Grassmannian.

We first argue that the number of solutions to the enumerative problem is
the degree $\deg{\mathcal K}^q_{m,p}$ of the quantum Grassmannian, and then
use the form of a Gr\"obner basis for the ideal 
${\mathcal I}^q_{m,p}={\mathcal I}({\mathcal K}^q_{m,p})$ of the quantum
Grassmannian to give another proof of 
Theorem~\ref{degree}, that $\deg{\mathcal K}^q_{m,p}$ is the number of
maximal chains in the poset ${\mathcal C}^q_{m,p}$ of quantum Pl\"ucker
coordinates.

\subsection{The enumerative problem of Question 1}

Given an $m$-plane $L\subset{\mathbb K}^{m+p}$ and a point
$s\in{\mathbb P}^1$, the set of degree $q$ maps
$M\in{\mathcal M}^q_{m,p}$ which satisfy
\begin{equation}\label{HypSec}
  M(s)\cap L\ \neq\ \{0\}
\end{equation}
is a hyperplane section of ${\mathcal M}^q_{m,p}$ in its Pl\"ucker
embedding.
This was shown both in Section 2.3 and, for special versions
of~(\ref{HypSec}), in Section 3.
Thus the enumerative problem of Question 1 asks for the number of points
common to ${\mathcal M}^q_{m,p}$ and to 
$N:=q(m+p)+mp (\,=\dim {\mathcal M}^q_{m,p})$ hyperplanes.
Hence if there are finitely many solutions, their number is bounded by 
$\deg{\mathcal K}^q_{m,p}$ and it equals this degree if there are no
points in the boundary ${\mathcal K}^q_{m,p}-{\mathcal M}^q_{m,p}$
common to all the hyperplanes.

Given a point $s\in{\mathbb P}^1$, the {\it evaluation map}
$\mbox{\rm ev}_s\colon{\mathcal M}^q_{m,p}\to
 \mbox{\it Grass}(p,{\mathbb K}^{m+p})$ associates a map $M$
to the $p$-plane $M(s)$.
The evaluation map extends to the quantum Grassmannian.
One way to see this is that the evaluation map is defined on the 
Quot scheme ${\mathcal Q}^q_{m,p}$ and it factors through the Pl\"ucker map
${\mathcal Q}^q_{m,p}\to{\mathcal K}^q_{m,p}$~\cite{Be97}.
Concretely, points $M$ of ${\mathcal K}^q_{m,p}$ are represented (possibly
non-uniquely) by matrices of homogeneous forms whose minors are forms of
degree $q$, and a point is on the boundary if the minors are not relatively
prime.
The (classical) Pl\"ucker coordinates of 
$\mbox{\rm ev}_s(M)\in\mbox{\it Grass}(p,{\mathbb K}^{m+p})$ are 
given by first dividing each minor by the common polynomial factor, and then
evaluating at the point $s\in{\mathbb P}^1$.

Since $\mbox{\it Grass}(p,{\mathbb K}^{m+p})$ is a homogeneous space,
Kleiman's Properness Theorem~\cite[Theorem 2(i)]{MR50:13063} shows that 
for distinct points $s_1,s_2,\ldots,s_N\in{\mathbb P}^1$
and general $m$-planes
$L_1,L_2,\ldots,L_N\in{\mathbb K}^{m+p}$, the collection of
hyperplanes given by
$$
   M(s_i)\cap L_i\ \neq\ \{0\}
$$
meet the boundary properly.
Since the boundary has dimension $N-m-p+1$, there are no points common to
the boundary and these hyperplanes.

Again by Kleiman's Properness Theorem, there are finitely many points $M$ in 
${\mathcal K}^q_{m,p}$ (and hence in ${\mathcal M}^q_{m,p}$) common to these
hyperplanes, as $N$ is the dimension of ${\mathcal K}^q_{m,p}$.
Since this set is a particular complementary linear section of 
${\mathcal K}^q_{m,p}$, its number of points
(possibly counted with multiplicity) is $\deg{\mathcal K}^q_{m,p}$.
When ${\mathbb K}$ is algebraically closed of characteristic zero, 
Kleiman's Transversality
Theorem~\cite[Theorem 2(ii)]{MR50:13063} implies that the solutions appear
without multiplicity, and so $\deg{\mathcal K}^q_{m,p}$ solves the
enumerative problem of Question 1.

\subsection{The degree of ${\mathcal K}^q_{m,p}$ via Gr\"obner bases}

For basics on Gr\"obner bases, we recommend either~\cite{CLO92}
or~\cite{Sturmfels_GBCP}, whose Chapter 11 has a description of the
classical $q=0$ version of the results discussed here.
Let ${\mathbb K}[{\mathcal C}^q_{m,p}]$ be the ring generated by the quantum
Pl\"ucker coordinates $z_{\alpha^{(a)}}$ for 
$\alpha^{(a)}\in{\mathcal C}^q_{m,p}$, the coordinate ring of the
Pl\"ucker space 
${\mathbb P}(\bigwedge^p{\mathbb K}^{m+p}\otimes{\mathbb K}^{q+1})$.
Let $\prec$ be the degree reverse lexicographic term order on 
ring ${\mathbb K}[{\mathcal C}^q_{m,p}]$ 
induced by an ordering of the variables
$z_{\alpha^{(a)}}$ corresponding to any (fixed) linear extension of the
poset ${\mathcal C}^q_{m,p}$.
The poset ${\mathcal C}^q_{m,p}$ is in fact a distributive lattice,
with $\alpha^{(a)}\wedge\beta^{(b)}$ the meet (greatest lower bound) 
and $\alpha^{(a)}\vee\beta^{(b)}$ the join (least upper bound) of the
indices $\alpha^{(a)}$ and $\beta^{(b)}$.

\begin{theorem}[\cite{SS_SAGBI}]\label{thm:gbasis}
The reduced Gr\"obner basis of the Pl\"ucker ideal 
${\mathcal I}^q_{m,p}$ of the quantum Grassmannian ${\mathcal K}^q_{m,p}$
consists of quadratic polynomials in 
$\,{\mathbb K}[{\mathcal C}^q_{p,m}]\,$ 
 which are indexed by 
pairs of incomparable variables $\gamma^{(c)},\delta^{(d)}$ in the poset 
$\,{\mathcal C}^q_{p,m}$,
$$
 S(\gamma^{(c)},\delta^{(d)}) \quad = \quad
  z_{\gamma^{(c)}}\cdot z_{\delta^{(d)}}\ -\ 
  z_{\gamma^{(c)}\vee\delta^{(d)}}\cdot 
  z_{\gamma^{(c)}\wedge\delta^{(d)}}
\,\, + \,\,\hbox{lower terms in $\prec$},
$$
and all lower terms $\,\lambda z_{\alpha^{(a)}}z_{\beta^{(b)}}\,$ 
of $\,S(\gamma^{(c)},\delta^{(d)}) \,$
satisfy $\,\alpha^{(a)}<\gamma^{(c)}\wedge\delta^{(d)}$ and 
$\gamma^{(c)}\vee\delta^{(d)}<\beta^{(b)}$.
\end{theorem}

By Theorem~\ref{thm:gbasis}, the initial ideal
$\mbox{\rm in}_\prec({\mathcal I}^q_{m,p})$ of the Pl\"ucker ideal
is generated by all monomials $z_{\alpha^{(a)}}z_{\beta^{(b)}}$ 
with $\alpha^{(a)},\beta^{(b)}\in{\mathcal C}^q_{m,p}$ incomparable.
We write this initial ideal as an intersection of prime 
ideals.
For this, let $\mbox{\CYR CH}^q_{m,p}$
be the set of (saturated) chains in the poset ${\mathcal C}^q_{m,p}$.

\begin{lemma}\label{lem:pdp}
Let ${\mathcal I}^q_{m,p}$ be the Pl\"ucker ideal.
Then
$$
  \mbox{\rm in}_\prec({\mathcal I}^q_{m,p})\ =\ 
  \bigcap_{\mbox{\Cyr ch}\in \mbox{\Cyr CH}^q_{m,p}C^q_{m,p}}
  \langle  z_{\delta^{(d)}}\mid \delta^{(d)}\not\in \mbox{\CYR ch}\rangle\,.
$$
\end{lemma}

\begin{proof}
If $z_{\alpha^{(a)}} z_{\beta^{(b)}}$ is a generator of 
$\mbox{\rm in}_\prec({\mathcal I}^q_{m,p})$, then 
$\alpha^{(a)}$
and $\beta^{(b)}$ are incomparable in 
the poset ${\mathcal C}^q_{m,p}$.
Thus if $\mbox{\CYR ch}\in \mbox{\CYR CH}^q_{m,p}$ is a saturated chain, at most
one of $\alpha^{(a)}$ or $\beta^{(b)}$ lies in the chain $\mbox{\CYR ch}$, and so 
$z_{\alpha^{(a)}} z_{\beta^{(b)}}$ lies in the ideal
$\langle z_{\delta^{(d)}}\mid \delta^{(d)}\not\in \mbox{\CYR ch}\rangle$.

Suppose now that $z$ is a monomial not in the initial ideal
$\mbox{\rm in}_\prec({\mathcal I}^q_{m,p})$.
Then the variables appearing in $z$ have indices which are comparable
in the poset ${\mathcal C}^q_{m,p}$.
Thus we may write 
$z=z_{\alpha^{(a)}}\cdot z_{\beta^{(b)}}\cdots z_{\gamma^{(c)}}$
with $\alpha^{(a)}\leq \beta^{(b)}\leq\cdots\leq\gamma^{(c)}$.
There is some chain $\mbox{\CYR ch}\in\mbox{\CYR CH}^q_{m,p}$ containing the
indices $\alpha^{(a)},\beta^{(b)},\ldots,\gamma^{(c)}$ and so the monomial 
$z$ does not lie in
the ideal $\langle z_{\delta^{(d)}}\mid \delta^{(d)}\not\in\mbox{\CYR ch}\rangle$.
This proves the equality of the two monomial ideals.
\end{proof}

Each ideal $\langle z_{\delta^{(d)}}\mid \delta^{(d)}\not\in\mbox{\CYR ch}\rangle$
defines the coordinate subspace of Pl\"ucker space spanned by the
coordinates $z_{\alpha^{(a)}}$ with $\alpha^{(a)}\in\mbox{\CYR ch}$, which is
isomorphic to ${\mathbb P}^{q(m+p)+mp}$, as every maximal chain $\mbox{\CYR ch}$
of  ${\mathcal C}^q_{m,p}$ has length 
$q(m+p)+mp+1$.
Thus the zero scheme of 
$\mbox{\rm in}_\prec({\mathcal I}^q_{m,p})$ is the union of these
coordinate subspaces, and so it has degree equal to their number.

\begin{proof}[Alternative Proof of Theorem~\ref{degree}]
The degree of the quantum Grassmannian is the degree of its 
ideal ${\mathcal I}^q_{m,p}$.
By Macaulay's Theorem~\cite{Mac1927} (see also~\cite[\S1.10]{E95}), this is
the degree of the initial  
ideal, $\deg(\mbox{\rm in}_\prec({\mathcal I}^q_{m,p}))$, which is equal
to the number of chains in $\mbox{\CYR CH}^q_{m,p}$, by Lemma~\ref{lem:pdp}.
\end{proof}

\begin{remarks}\mbox{ }
\begin{enumerate}
\item[(1)]
In~\cite{SS_SAGBI}, reduced Gr\"obner bases for the quantum Schubert
varieties $Z_{\alpha^{(a)}}$ which are restrictions of the Gr\"obner basis of
Theorem~\ref{thm:gbasis} are also constructed, and a consequence is that the
definition~(\ref{qsv}) is in fact ideal-theoretic:
$$
  {\mathcal I}(Z_{\alpha^{(a)}})\ =\ {\mathcal I}^q_{m,p} + 
    \langle z_{\beta^{(b)}}\mid\beta^{(b)}\not\leq \alpha^{(a)}\rangle\,.
$$
The form of these Gr\"obner bases also significantly strengthens
Proposition~\ref{prop:RRW} (ii)  to the level of homogeneous ideals.

\item[(2)]
The reduced Gr\"obner basis for the Pl\"ucker ideal of the classical
Grassmannian $(q=0)$ may be constructed as follows~\cite{HP52,Sturmfels_invariant}:
First a Gr\"obner basis consisting of linearly independent quadratic
polynomials, one for each incomparable pair, is constructed using invariant
theory. 
Then this basis is reduced to obtain the desired reduced Gr\"obner basis.
In contrast to that approach, the reduced Gr\"obner basis of
Theorem~\ref{thm:gbasis} was constructed explicitly using a double induction
on the poset ${\mathcal C}^q_{m,p}$.

For $\alpha^{(a)}\geq\beta^{(b)}$ in ${\mathcal C}^q_{m,p}$,
define the skew quantum Schubert variety
$$
 \qquad Z_{\alpha^{(a)}/\beta^{(b)}}\quad :=\quad 
  \{(z_{\gamma^{(c)}})\in {\mathcal K}^q_{m,p}\colon
  z_{\gamma^{(c)}}=0\ \mbox{ unless }\ 
     \alpha^{(a)}\geq \gamma^{(c)}\geq \beta^{(b)} \}\ .
$$
A first step is to show 
that if $\alpha^{(a)}=\gamma^{(c)}\vee\delta^{(d)}$, then 
$$
  z_{\gamma^{(c)}}\cdot z_{\delta^{(d)}}\ -\ 
  z_{\gamma^{(c)}\vee\delta^{(d)}}\cdot 
  z_{\gamma^{(c)}\wedge\delta^{(d)}}\ \in\ 
   {\mathcal I}(Z_{\alpha^{(a)}/\beta^{(b)}})\,,
$$
by downward induction on $\beta^{(b)}$.
When $\beta^{(b)}\in {\mathcal C}^q_{m,p}$ is minimal,
$Z_{\alpha^{(a)}/\beta^{(b)}}=Z_{\alpha^{(a)}}$.
Then the forms $S(\gamma^{(c)},\delta^{(d)})$ of Theorem~\ref{thm:gbasis}
are constructed by increasing induction on $\alpha^{(a)}$.

\item[(3)]  
An important part of~\cite{SS_SAGBI} was to study the rational
parameterization of  
${\mathcal K}^{pn}_{m,p}$ given by $m{+}p$ by $p$ matrices whose entries are
generic polynomials of degree $n$, and also by an intermediate variety,
the Grassmannian of $p$-planes in ${\mathbb K}^{(pn+1)(m+p)}$.
This `long Grassmannian' was used by Byrnes~\cite{By83} to obtain a
different compactification of ${\mathcal M}^q_{m,p}$ than 
${\mathcal K}^q_{m,p}$.
It was also used to prove Proposition~\ref{prop:RRW}(ii)~\cite{RRW96,RRW98},
and the indices 
of its Schubert varieties appeared implicitly in the indexing scheme of
Section 2.5.
Lastly, the classical (ideal-theoretic) version of  
Proposition~\ref{prop:RRW}(ii) for Schubert varieties in the long Grassmannian
was used in the inductive steps of item (2) above.

\item[(4)]
We expect this approach and these results to generalize to other flag manifolds,
giving an analog of standard monomial theory~\cite{LaSe91} for spaces of
rational curves in all flag manifolds.
\end{enumerate}

\end{remarks}


\section{Quantum Cohomology and the Formula of Vafa and Intriligator}\label{sec:last}

We describe some of the standard story of the enumerative problem of
Question~1.
We first briefly review some history of the formula of Vafa and
Intriligator.
Next, we visit the classical cohomology ring of the Grassmannian and its
quantum deformation, and then give the formula of Vafa and
Intriligator.
We then show how this formula of  Vafa and Intriligator agrees with the
formula~(\ref{formula}) of Ravi, Rosenthal, and Wang.
We next give an alternative way to view the quantum cohomology ring of the
Grassmannian, and discusses how this same ring arose in two different contexts in
representation theory. 
This survey concludes with some open problems concerning quantum Littlewood-Richardson
coefficients. 
\smallskip

Inspired by Donaldson's invariants of 4-manifold~\cite{Do87}, 
Gromov~\cite{Gr85} proposed that topological invariants of moduli spaces of
pseudo-holomorphic curves in a symplectic manifold $X$ would give invariants of
the symplectic structure of $X$.
Following ideas of Witten~\cite{Wi91}, Vafa~\cite{Va92} proposed so-called quantum
multiplications in the cohomology rings of symplectic manifolds with structure
constants certain correlation functions, and
conjectured remarkable residue formulae for these correlation functions when $X$ is a
Grassmannian. 
This was made more precise by Intriligator~\cite{In91}.
Ruan (see~\cite{Ru98}) was perhaps the first to link this work in theoretical
physics to the work of Gromov, realizing that Witten's correlation functions
were in fact Gromov's invariants, and hence the formula of Vafa and
Intriligator computes intersection numbers of curves of all genera on
Grassmannians.
Siebert and Tian~\cite{ST97} generalized the program of Vafa and Intriligator
from the Grassmannian 
to certain Fano manifolds---in particular, they proved the formula of Vafa
and Intriligator and constructed the (small) quantum cohomology rings of these
manifolds.
Previously (and with different methods), Bertram,  Daskalopoulos, and
Wentworth~\cite{BDW96} had proven this formula for genus 1 invariants of
high degree curves in Grassmannians of 2-planes, and Bertram~\cite{Be97}
later developed a quantum Schubert calculus which enabled the computation of
intersection numbers involving arbitrary Schubert conditions.

\subsection{The cohomology ring of the Grassmannian}
The cohomology ring of the complex Grassmannian 
$\mbox{\it Grass}(p,{\mathbb C}^{m+p})$ has a standard
presentation
 \begin{equation}\label{BT}
  H^*(\mbox{\it Grass}(p,{\mathbb C}^{m+p}))
   \ \stackrel{\sim}{\longrightarrow}\ 
  {\mathbb C}[c_1,c_2,\ldots,c_p]/
   \langle h_{m+1},h_{m+2},\ldots,h_{m+p} \rangle\,,
 \end{equation}
where $\deg c_i=2i$ and 
$h_1,h_2,\ldots,h_{m+p}$ are defined recursively in terms of the
$c_i$ as follows
 \begin{equation}\label{eq:Hs}
   h_j-c_1h_{j-1}+\cdots+(-1)^{j-1}c_{j-1}h_1+(-1)^jc_j\ =\ 0\,,
 \end{equation}
with $c_i=0$ for $i>p$.
The isomorphism is given by associating $c_i$ to the $i$th Chern class of the dual $S^*$ of
the tautological rank $p$ subbundle $S$ over the Grassmannian.
Then $h_i$ is the $i$th Chern class of the rank $m$ quotient
bundle $Q$, and these classes vanish for $i>m$.
The relation~(\ref{eq:Hs}) between these classes $c_i$ and $h_i$ is succinctly expressed
via the splitting principle,
$$
  1\ =\ c({\mathbb C}^{m+p})\ =\ c(S^*) c(Q^*)\,,
$$
where ${\mathbb C}^{m+p}$ is the trivial bundle and $c(\cdot)$ is the total Chern class.
(Here, $c(Q^*)=1-h_1+h_2-\cdots$.)

Because the cohomology ring is a complete intersection and the $h_j$ are
homogeneous of degree $2j$, it is Gorenstein with socle in dimension
$2mp=\sum_j(\deg h_{m+j}-\deg c_j)$.
A generator of the socle is the image of $c^m_p$, and the degree map 
(used to compute intersection numbers) is simply the coefficient of 
$c^m_p$ in an element of this quotient ring.
Thus, given some classes $\xi_1,\xi_2,\ldots,\xi_l$ in cohomology which are 
Poincar\'e
dual to cycles $X_1,X_2,\ldots,X_l$ in general position, the
coefficient of $c^m_p$ in the product 
$\xi_1\cdot\xi_2\cdots\xi_l$ is the number of points in the intersection
$$
  X_1\cap X_2\cap\cdots\cap X_l\,,
$$
when there are finitely many such points.

What is less known is that the degree map may be computed using the local
residue associated to the map 
${\bf H}:=((-1)^{p-1}h_{m+1},(-1)^{p-2}h_{m+2},\ldots,h_{m+p})\colon
     {\mathbb C}^p\to{\mathbb C}^p$.
(See~\cite[\S4]{ST97} for details.)
This residue is
$$
  \mbox{\rm Res}_{\bf H}(F)\ =\ 
   \frac{1}{(2\pi i)^p}\int_{\Gamma^\epsilon}
    \frac{F}{(-1)^{p-1}h_{m+1}(-1)^{p-2}h_{m+2}\cdots h_{m+p}}dc_1\cdot dc_2\cdots dc_p\,,
$$
for $F\in{\mathbb C}[c_1,c_2,\ldots,c_p]$.
Here $\Gamma^\epsilon$ is a smooth canonically oriented cycle in the 
region where no component $h_{m+i}$ of ${\bf H}$ vanishes.
Standard properties of residues~\cite[\S5]{GH78} imply that the residue
vanishes on the ideal of~(\ref{BT}), and so gives a well-defined map on the
cohomology ring.
Furthermore, when $F$ is homogeneous, the residue vanishes unless 
$\deg F=2mp$, for otherwise the form is exact.
Thus the residue is proportional to the degree map, and the calculations
we do below show the constant of proportionality is $(-1)^{\binom{p}{2}}$.

The presentation~(\ref{BT}) has another form.
Let $W=\frac{1}{m+p+1}P_{m+p+1}$, where $P_{m+p+1}$ is the 
$(m{+}p{+}1)$th Newton power sum symmetric polynomial in the variables 
$x_1,\ldots,x_p$.
If we express $P_{m+p+1}$ as a polynomial in the elementary symmetric
polynomials $c_1,c_2,\ldots,c_p$, then we have (see below)
 \begin{equation}\label{eq:power_sum}
  (-1)^{1+j}h_{m+p+1-j}\ =\ \frac{\partial W}{\partial c_j}\,,
 \end{equation}
where $h_i$ is the $i$th complete homogeneous symmetric polynomial
in the variables $x_j$ (these satisfy~(\ref{eq:Hs}) when the $c_i$ are elementary
symmetric polynomials).
Thus the presentation becomes
$$
  H^*(\mbox{\it Grass}(p,{\mathbb C}^{m+p}))
   \ \stackrel{\sim}{\longrightarrow}\ 
  {\mathbb C}[c_1,c_2,\ldots,c_p]/
   \langle dW=0 \rangle\,.
$$\smallskip

We derive~(\ref{eq:power_sum}), working in the ring $\Lambda$ of symmetric
functions in the indeterminates $x_1,x_2,\ldots$~\cite[Sect.~7]{St_ECII}~\cite[I.2]{Mac}~\cite{Sa91}. 
To obtain this formula for polynomials,  specialize $x_i$ to $0$ for $i>p$.
First, we have the fundamental relations
 \begin{eqnarray}
 1&=& \left(\sum_{r\geq 0} h_r (-t)^r\right) \cdot
       \left(\sum_{r\geq 0} e_r t^r\right)\,,\label{eq:E-H}\\
 \sum_{r>0} p_r (-t)^{r-1}&=& \frac{d}{dt} \log 
    \left(\sum_{r\geq 0} e_r t^r\right)\,,\label{eq:dlogE}
 \end{eqnarray}
where $e_i,h_i$, and $p_i$ are, respectively,  the elementary, complete homogeneous, and power
sum symmetric functions of degree $i$.
(Note that~(\ref{eq:E-H}) gives~(\ref{eq:Hs}).)
Differentiating~(\ref{eq:dlogE}) with respect to $e_j$ gives
$$
  \sum_{r>0} \frac{\partial p_r}{\partial e_j} (-t)^{r-1} \ =\ 
  \frac{d}{dt} \frac{t^j}{\sum_{r\geq 0} e_r t^r}\ =\ 
  \frac{d}{dt} t^j\sum_{r\geq 0} h_r (-t)^r\,.
$$
Equating coefficients of $t^{m+p}$ gives
$$
  (-1)^{m+p}\frac{\partial p_{m+p+1}}{\partial e_j}\ =\ 
  (-1)^{m+p+1-j} (m+p+1)h_{m+p+1-j}\,,
$$
from which~(\ref{eq:power_sum}) follows.

\subsection{Quantum cohomology and the formula of Vafa and
   Intriligator}\label{sec:VI} 
The quantum cohomology ring is a perturbation (depending on a K\"ahler form)
of the classical cohomology ring whose structure encodes the genus zero
Gromov-Witten invariants. 
For the Grassmannian, Vafa and Intriligator began with the
perturbation of $W$
$$
  QW\ :=\ W + (-1)^p\beta c_1\,,
$$
where $\beta$ is a complex number associated to the perturbing K\"ahler
form. 
For our enumerative problem, $\beta=1$.
They then proposed the following presentation
for the quantum cohomology ring 
$$
   QH^*(\mbox{\it Grass}(p,{\mathbb C}^{m+p}))
   \ =\   {\mathbb C}[c_1,c_2,\ldots,c_p]/ \langle dQW=0 \rangle
$$
which is 
 \begin{equation}\label{eq:qcoh-pres}
    {\mathbb C}[c_1,c_2,\ldots,c_p]/
       \langle (-1)^{p-1}h_{m+1}, (-1)^{p-2} h_{m+2}, \ldots,
            h_{m+p} + (-1)^p \rangle\,.
 \end{equation}

They also proposed the following formula.
Let $X_1,X_2,\dots,X_l$ be special Schubert cycles in the Grassmannian 
which are in general position.
Suppose $\beta=1$.
For a genus $g\geq 0$, set $\langle X_1,X_2,\ldots,X_l\rangle_g:=0$
unless the sum of the codimensions of the $X_j$ 
is equal to $d(m+p)+mp(1-g)$, for some non-negative integer $d$.
When there is such an integer $d$, let the Gromov-Witten invariant 
$\langle X_1,X_2,\ldots,X_l\rangle_g$ be the number of maps
$$
  f\ \colon\ (\Sigma,s_1,s_2,\ldots,s_l)\ \longrightarrow\ 
    \mbox{\it Grass}(p,{\mathbb C}^{m+p})
$$
satisfying $f(s_i)\in X_i$.
Here $\Sigma$ is a fixed genus $g$ curve, $s_1,s_2,\ldots,s_l$ are fixed, but
general, points of  $\Sigma$, and $f_*[\Sigma]=d\cdot c_1$.
Determining when this definition is well-founded and providing a satisfactory
alternative when it is not is an important and subtle story which we do not
relate. 
When $\beta\neq 1$, the definition involves
pseudo-holomorphic curves, and we omit it.

Suppose we have special Schubert classes
$c_{i_1},c_{i_2},\ldots,c_{i_l}$ with $c_{i_j}$ Poincar\'e dual to 
$X_j$.
Then the formula of Vafa and Intriligator for 
$\langle X_1,X_2,\ldots,X_l\rangle_{g}$ is
 \begin{equation}\label{VI}
    (-1)^{\binom{p}{2}(g-1)}\sum_{dQW(c_1,c_2,\ldots,c_p)=0}
       \det\left(\frac{\partial^2QW}{\partial c_i\partial c_j}\right)^{g-1}
        \cdot c_{i_1}\cdot c_{i_2}\cdots c_{i_l}\,.
 \end{equation}
One remarkable feature of this formula involves the determinant 
$J=\det\left(\frac{\partial^2QW}{\partial c_i\partial c_j}\right)$.
The formula implies that the genus $g$ Gromov-Witten invariant of a monomial 
$c^{\bf i}$ equals the genus $g-1$ Gromov-Witten invariant of
$c^{\bf i}/J$ (up to a sign).

We relate this to the classical intersection formula when $g=0$.
Since $dQW$ is the vector
$((-1)^{p-1}h_{m+1},\ldots,-h_{m+p-1},h_{m+p}+(-1)^p\beta)$,
the determinant $J$ in the formula~(\ref{VI}) is also the Jacobian of the 
map ${\bf H}$, and so the summand of~(\ref{VI}) becomes
$$
   \sum_{c\in {\bf H}^{-1}(y)}  \frac{F(c)}{J}\,,
$$
where  $y=(0,\ldots,0,(-1)^{p+1}\beta)$ and 
$F(c)=c_{i_1}\cdot c_{i_2}\cdots c_{i_l}$.
Let $\mbox{\rm res}_y(F)$ denote this number, which is a trace but also a
residue as $y$ is a regular value of the map ${\bf H}$.
A further property of the residue is that $\mbox{\rm res}_y(F)$ extends
holomorphically to a neighborhood of 0, and
$$
  \lim_{y\to 0}\mbox{\rm res}_y(F)\ =\ \mbox{\rm Res}_{\bf H}(F)\,.
$$
This shows rather explicitly how this formula of Vafa and Intriligator
is a deformation of the classical intersection formula.

\subsection{Relation between the formulae of Vafa and Intriligator and of 
  Ravi, Rosenthal, and Wang}

We relate the formula~(\ref{VI}) for genus 0 curves to the
formula~(\ref{formula}) of Ravi, Rosenthal, and Wang.
For $\alpha\in\binom{[n]}{p}$, let $\Omega_\alpha=Z_{\alpha^{(0)}}$, a
Schubert subvariety of $\mbox{\it Grass}(p,{\mathbb C}^{m+p})$.

\begin{theorem}[\cite{RRW98}]
 Let $\alpha^{(a)}\in{\mathcal C}^q_{m,p}$.
 Then
$$
  \deg Z_{\alpha^{(a)}}\ =\ (-1)^{\binom{p}{2}} \sum_{dQW=0} 
      \frac{\sigma_{\alpha^\vee}\cdot c_1^{|\alpha^{(a)}|}}{J}\,,
$$
 where $\sigma_{\alpha^\vee}$ is the cohomology class Poincar\'e dual to the
 fundamental cycle of $\Omega_\alpha$.
\end{theorem}

\noindent{\it Proof.}
Since the Gromov-Witten invariants of genus zero curves on the Grassmannian
may be computed in the 
quantum cohomology ring of the Grassmannian, the obvious linear extension of
the formula of Vafa and Intriligator~(\ref{VI}) for genus zero curves to
arbitrary cycles $X_i$ is valid.
Thus the right hand side above computes the Gromov-Witten invariant
 \begin{equation}\label{G-W_degree}
  \langle \Omega_\alpha,\, X_1,\, X_2,\, \ldots,
    \,X_{|\alpha^{(a)}|}\rangle_0\,,
 \end{equation}
where $X_1,X_2,\ldots,X_{|\alpha^{(a)}|}$ are special Schubert varieties in
general position, each dual to $c_1$.
Since the cohomological degree of the class 
$\sigma_{\alpha^\vee}{\cdot}c_1^{|\alpha^{(a)}|}$ is 
 \begin{eqnarray*}
  2(mp-|\alpha|)+2|\alpha^{(a)}| &=& 2mp-2|\alpha|+2(a(m+p)+|\alpha|)\\
          &=& 2a(m+p) + 2mp\,,
 \end{eqnarray*}
this is an invariant of degree $a$ curves.

We first express this Gromov-Witten invariant as the number of points in 
an intersection.
Given a point $s\in{\mathbb P}^1$, the evaluation map
ev$_s\colon{\mathcal M}^a_{m,p}\to\mbox{\it Grass}(P,{\mathbb C}^{m+p})$
associates a curve $M$ to the $p$-plane $M(s)$.
Each cycle $X_i$ has the form
$$
  \Omega_{L_i}\ =\ \{H\in{\it Grass}(p,{\mathbb C}^{m+p})\mid 
                      H\cap L_i\neq\{0\}\}\,,
$$
where $L_1,\ldots,L_{|\alpha^{(a)}|}$ are $m$-planes in 
${\mathbb C}^{m+p}$ in general position.
Thus~(\ref{G-W_degree}) counts the number of points in the
intersection
 \begin{equation}\label{q-int}
   \mbox{\rm ev}_{s_0}^{-1}\Omega_\alpha\bigcap
    \mbox{\rm ev}_{s_1}^{-1}\Omega_{L_1}\bigcap\cdots\bigcap
     \mbox{\rm ev}_{s_{|\alpha^{(a)}|}}^{-1}\Omega_{L_{|\alpha^{(a)}|}}\,,
 \end{equation}
where $s_0,s_1,\ldots,s_{|\alpha^{(a)}|}$ are general points in 
${\mathbb P}^1$.

Observe that we may choose $s_0$ to be the point $\infty$ at infinity in  
${\mathbb P}^1$.
Then, in the quantum Pl\"ucker coordinates 
$(z_{\beta^{(b)}}\mid \beta^{(b)}\in{\mathcal C}^a_{m,p})$ for a point
$z\in {\mathcal M}^a_{m,p}$, 
$$
  \mbox{\rm ev}_\infty\ \colon\  
    z \ \longmapsto\ 
   (z_{\beta^{(a)}}\mid \beta\in{\textstyle \binom{[m+p]}{p})}\,.
$$
Since, in the Pl\"ucker coordinates $(y_\beta\mid\beta\in\binom{[m+p]}{p}))$
for $\mbox{\it Grass}(p,{\mathbb C}^{m+p})$ we have the 
analog of~(\ref{qsv}) for $\Omega_\alpha$,
$$
  \Omega_\alpha\ :=\ 
    \{y\in\mbox{\it Grass}(p,{\mathbb C}^{m+p}\mid 
     y_\beta=0\mbox{ if }\beta\not\leq\alpha\}\,,
$$
we see that 
$\overline{{\rm ev}_\infty^{-1}(\Omega_\alpha)}=Z_{\alpha^{(a)}}$.

Finally, observe that for an $m$-plane $L\subset{\mathbb C}^{m+p}$
and a point $s\in{\mathbb P}^1$, 
$$
   {\rm ev}_s^{-1}\Omega_L\ =\ 
    \{M\in{\mathcal M}^a_{m,p}\mid M(s)\cap L\neq\{0\}\}\,
$$
which is a hyperplane section of ${\mathcal M}^a_{m,p}$.
Thus the number of points in the intersection~(\ref{q-int}) is bounded by the
degree of $Z_{\alpha^{(a)}}$ and it equals this degree if all points of
intersection lie in 
$Z_{\alpha^{(a)}}\cap{\mathcal M}^a_{m,p}$.
But this occurs for general $m$-planes $L_i$ and points $s_i$,
by Remark~\ref{rem:intersection}.
\qed\medskip

This proof is unsatisfactory in that both
sides of the equation have a simple algebraic-combinatorial interpretation, yet
we argued using the definition of the Gromov-Witten invariants, rather
than something more elementary.
We now give a more direct proof, following~\cite{RRW98}.
For a sequence $I:0<i_1<i_2<\cdots<i_p$ of integers, let $S_I$ be the Schur symmetric
polynomial in $x_1,\ldots,x_p$ associated to the partition
$(i_p-p,\ldots,i_2-2, i_1-1)$, which is also a polynomial in the 
elementary symmetric polynomials.

\begin{theorem}\label{thm:last}
 For $\alpha^{(a)}\in{\mathcal C}^q_{m,p}$, define
$$
  \delta(\alpha^{(a)})\ :=\  (-1)^{\binom{p}{2}} \sum_{dQW=0}
    \frac{ c_1^{|\alpha^{(a)}|} \cdot S_{\alpha^\vee}}{J}\,.
$$
 Then the function $\delta(\alpha^{(a)})$ satisfies the
 recursion~(\ref{recursion}). 
\end{theorem}

This will prove the equality of the two formulae, since
under the map~(\ref{BT}) to cohomology we have
$$
  S_\alpha\ \longmapsto \sigma_\alpha
   \qquad\mbox{\rm for }\alpha\in{\textstyle \binom{[m+p]}{p}}\,.
$$

\noindent{\it Proof.}
We set $n:=m+p$ and change coordinates, working in the ring
of symmetric polynomials ${\mathbb C}[x_1,x_2,\ldots,x_p]^{{\mathcal S}_p}$ 
in $x_1,x_2,\ldots,x_p$.
If we let each $x_i$ have cohomological degree 2, then this ring is isomorphic
to the ring 
${\mathbb C}[c_1,c_2,\ldots,c_p]$ with the isomorphism given by
$c_i=e_i(x_1,x_2,\ldots,x_p)$.
Here $e_i(x_1,x_2,\ldots,x_p)$ is the $i$th elementary symmetric polynomial in 
$x_1,x_2,\ldots,x_p$.

This theorem is a consequence of Lemma~\ref{lem:recursion} and the following
lemma.
Let $y_1,y_2,\ldots,y_n$ be the $n$th roots of $(-1)^{p+1}$.

\begin{lemma}\label{lem:funD}
For $K=k_1,k_2,\ldots,k_p$, define the function $D(K)$ to be 
 \begin{equation}\label{fnD}
   \frac{(-1)^{\binom{p}{2}}}{n^p}\sum
   (x_1\cdots x_p)(x_1+\cdots+x_p)^{\sum_j k_j-j}
   \det(x_j^{n-k_i})\,\det(x_i^{p-j})\,,
 \end{equation}
the sum over all $I\in\binom{[n]}{p}$, where $x_j=y_{i_j}$.
Then
\begin{enumerate}
 \item[(i)] For sequences of integers 
      $k_1\leq k_2\leq \cdots \leq k_p$ with $\sum_j k_j-j\geq 0$
      and $k_p<k_1+n$, the function $D(k_1,\ldots,k_p)$ satisfies the
      recursion, initial condition, and boundary conditions of
      Lemma~\ref{lem:recursion}. 
      In particular, $D(J(\alpha^{(a)}))=\deg Z_{\alpha^{(a)}}$.
 \item[(ii)]
      For $\alpha^{(a)}\in{\mathcal C}^q_{m,p}$, 
      $\delta(\alpha^{(a)})=D(J(\alpha^{(a)}))$. 
\end{enumerate}
\end{lemma}

\noindent{\it Proof of Lemma~\ref{lem:funD}(i).}
For a sequence $K=(k_1,k_2,\ldots,k_p)$ of integers, let 
$f(K):=\det(x_i^{n-k_j})$.
This determinant is the sum of terms
$f(\pi;K):={\rm sgn}(\pi)x_{\pi(1)}^{n-k_1}\cdots x_{\pi(p)}^{n-k_p}$
over all permutations $\pi$ of $\{1,2,\ldots,p\}$, where 
${\rm sgn}(\pi)$ is the sign of the permutation $\pi$.
Multiplying this term by $x_1+\cdots+x_p$ gives
$$
  (x_1+\cdots+x_p)\cdot f(\pi;K)\ =\ 
     \sum_{a=1}^p f(\pi;k_1,\ldots,k_a{-}1,\ldots,k_p)\,.
$$
Summing over all permutations $\pi$ gives the Pieri formula
$$
  (x_1+\cdots+x_p)\cdot f(K)\ =\ \sum_{a=1}^p
       f(k_1,\ldots,k_a-1,\ldots,k_p)\,,
$$
and thus $D(K)$ satisfies the recursion~(\ref{eq:rec}) of
Lemma~\ref{lem:recursion}.
Since $D(K)$ is antisymmetric in its arguments, it satisfies
the boundary condition~(\ref{A}).
If $k_p=k_1+n$, then the first and last rows of the matrix
$(x_j^{n-k_i})$ are the scalar multiples $(-1)^{p+1}$ of each other, and so
the function $D$ satisfies the boundary condition~(\ref{C}).

To show that $D$ satisfies the initial condition~(\ref{initial}) and the
remaining boundary condition~(\ref{B}), consider the values of
$D(k_1,\ldots,k_p)$ when $k_1<\cdots<k_p$, $0=\sum_j(k_j-j)$, and
$k_p<k_1+n$.
For such sequences $K$, we show that
 \begin{equation}\label{in-cond}
  D(K)\ =\ \left\{\begin{array}{lcl}
      1&\ &\mbox{\rm if } K=(1,2,\ldots,p)\,,\\
      0&  &\mbox{\rm otherwise}\,.\end{array}\right.
 \end{equation}
The first case of this is the initial condition~(\ref{initial}).
We deduce the boundary condition~(\ref{B}) from the second case
of~(\ref{in-cond}).

Let $J=j_1<j_2<\cdots<j_p$ be sequence of integers 
satisfying $\sum_ij_i-i\geq 0$ with $j_1=0$ and $j_p<n=j_1+n$.
Applying the recursion~(\ref{eq:rec}) $\sum_i(j_i-i)$ times to $D(J)$, and
the boundary conditions~(\ref{A}) and~(\ref{C}) shows that $D(J)$ is a sum
of terms $D(K)$ for $K$ satisfying the conditions for~(\ref{in-cond}), but
with $k_1\leq 0$.   
Thus every such term vanishes, and so $D(J)=0$.

We prove~(\ref{in-cond}), which will complete the proof of
Lemma~\ref{lem:funD}~(i).
For the sequences $K$ of~(\ref{in-cond}), we have
 \begin{eqnarray*}
  D(K)&=&\frac{(-1)^{\binom{p}{2}}}{n^p}
          \sum (x_1\cdots x_p)\det(x_j^{n-{k_i}})\det(x_i^{p-j})\\
      &=&\frac{(-1)^{\binom{p}{2}}}{n^p}
         \sum \det(x_j^{n-{k_i}})\det(x_i^{p+1-j})\,,
 \end{eqnarray*}
where the sum is over all (ordered) $p$-tuples $(x_1,\ldots,x_p)$ of the
$n$th roots $(y_1,\ldots,y_n)$ of $(-1)^{p+1}$.
We apply the Cauchy-Binet formula to this sum of products of
determinants to obtain
$$
  D(K)\ =\ \frac{(-1)^{\binom{p}{2}}}{n^p}
          \det\left(
            \left[\begin{array}{cccc}
             y_1^{n-k_1}&y_2^{n-k_1}&\cdots&y_n^{n-k_1}\\
               \vdots & \vdots & \ddots & \vdots\\
             y_1^{n-k_p}&y_2^{n-k_p}&\cdots&y_n^{n-k_p}
             \end{array}\right]
              \cdot
            \left[\begin{array}{ccc}
               y_1^p&\cdots&y_1\vspace{2pt}\\
               y_2^p&\cdots&y_2\\
              \vdots&\ddots&\vdots\\
               y_n^p&\cdots& y_n
             \end{array}\right] \right)\ .
$$
Expanding this product, we obtain
 \begin{equation}\label{power}
     D(K)\ =\ \frac{(-1)^{\binom{p}{2}}}{n^p}
          \det\left[\begin{array}{cccc}
              P(n+p-k_1)& \cdots& P(n+1-k_1)\\
               \vdots&\ddots&\vdots\\
              P(n+p-k_p)& \cdots& P(n+1-k_p)
              \end{array}\right] \ ,
 \end{equation}
where $P(b)$ is the sum of the $b$th powers of the $y_i$.
Since the $y_i$ are the $n$th roots of $(-1)^{p+1}$, we have
$$
  P(b)\ =\ \left\{\begin{array}{lcl}
              (-1)^{(p+1)a}\cdot n&\ &\mbox{\rm if } b=an\,,\\
                      0&  &\mbox{\rm otherwise}\,.
                  \end{array}\right.
$$

The $i$th row of the determinant~(\ref{power}) has at most one non-zero
entry, in the column $j$ where $j+k_i\equiv p+1$ modulo $n$.
Suppose that $K$ satisfies the conditions of~(\ref{in-cond}) and 
the determinant of~(\ref{power}) does not vanish.
Then each component of $K$ is congruent to one of $\{1,2,\ldots,p\}$ modulo 
$n$.
Since each congruence must occur for the determinant to be non-zero
(if you like, since no two components of $K$ are congruent modulo $n$),
we have that $K\equiv\{1,2,\ldots,p\}$ modulo $n$.
In particular, no component of $K$ vanishes.
Let $r$ be the index such that $k_r< 0< k_{r+1}$.
Since $k_p<k_1+n$ and $\sum_j(k_j-j)$, we must have $k_1\leq 1$, and also
$-n<k_1<k_p<n$.
In fact the condition that $K\equiv\{1,2,\ldots,p\}$ modulo $n$ implies that
$k_1,\ldots,k_r\leq -m$ and $0<k_{r+1},\ldots,k_p\leq p$. 
This implies $k_1\leq -m+1-r$ and $p{-}r\leq k_p$, and hence
$p-r<k_p<k_1+n\leq p+1-r$, which implies that
$K=(-m+1-r,-m+2-r,\ldots,-m,1,2,\ldots,p-r)$
and so $\sum_jk_j-j=-nr$.
Since this sum must equal 0, we see that 
$K$ must be $(1,2,\ldots,p)$.

When $K=(1,2,\ldots,p)$, the matrix of power sums is antidiagonal
with entries $(-1)^{p+1} n$, and so the determinant is
$(-1)^{p^2+p+\binom{p}{2}}n^p=(-1)^{\binom{p}{2}}n^p$, which shows that
$D(1,2,\ldots,p)=1$, as claimed.
This completes the proof of Lemma~\ref{lem:funD}(i).
\qed\medskip

\noindent{\it Proof of Lemma~\ref{lem:funD}(ii).}
We show that for $\alpha^{(a)}\in{\mathcal C}^q_{m,p}$, 
we have the equality $\delta(\alpha^{(a)})=D(J(\alpha^{(a)}))$. 
Since 
$$
  QW(x_1,\ldots,x_p)\ =\ 
     \sum_{j=1}^p \frac{x_j^{n+1}}{n+1} +(-1)^p x_j\,,
$$
the set of solutions for $dQW=0$ are just the set of $p$-tuples
of $n$th roots of $(-1)^{p+1}$.

The Schur polynomial $S_{\alpha^\vee}$ is equal to
the quotient of alternants~\cite{Ca1815}
$$
   \frac{\det(x_j^{n-\alpha_i})}{\det(x_i^{p-j})}\,.
$$
The denominator is the Vandermonde determinant
$\Delta:=\prod_{i<j}(x_i-x_j)$.

The Jacobian $J$ is the determinant of the Hessian of $QW$ with respect to the variables $c_i$,
which we compute using the multivariate chain rule
$$
  \left(\frac{\partial^2 QW}{\partial x_i\partial x_j}\right)\ =\ 
    \left(\frac{\partial^2 QW}{\partial c_i\partial c_j}\right)
    \cdot \left(\frac{\partial c_i}{\partial x_j}\right)^2\ +\ 
    \left( \sum_k \frac{\partial QW}{\partial c_k}
      \frac{\partial^2 c_k}{\partial x_i\partial x_j}\right)\ .
$$
Since we evaluate this where $dQW=0$, we obtain
$$
  \det\left.\left(\frac{\partial^2 QW}{\partial x_i\partial x_j}\right)\right|_{dQW=0}
    \ =\ 
  \det\left(\frac{\partial^2 QW}{\partial c_i\partial c_j}\right)
  \cdot\left[\det\left(\frac{\partial c_i}{\partial x_j}\right)\right]^2\,.
$$
The Hessian of $QW$ with respect to the variables $x_i$ is the diagonal matrix
with entry $nx_i^{n-1}$ in position $(i,i)$
and by Lemma~\ref{ToDo} below, 
$\det(\partial c_i/\partial x_j)=\Delta$.
Since  $\delta(\alpha^{(a)})$ is the sum over $p$-tuples 
$(x_1,\ldots,x_p)$ of $n$th roots of $(-1)^{p+1}$, we compute the 
value of the Jacobian $J=\det(\partial^2 QW/\partial c_i\partial c_j)$
at the $p$-tuple $(x_1,\ldots,x_p)$
to be
$$
  J\ =\ \frac{n^p(x_1\cdots x_p)^{n-1}}{\Delta^2}
   \ =\ \frac{n^p(x_1\cdots x_p)^n}{(x_1\cdots x_p)\Delta^2}
   \ =\ \frac{n^p}{(x_1\cdots x_p)\Delta^2}\,,
$$
as $(x_1\cdots x_p)^n=(-1)^{p(p+1)}=1$.

Since each summand involves the Vandermonde, we may restrict the sum
to be over the set $I$ of 
all $p$-tuples of distinct roots, which we will always take to be
in an order compatible with a fixed ordering of the $n$th roots
$y_1,\ldots,y_n$ of $(-1)^{p+1}$.
We may put these calculations together and obtain the following formula 
for $\delta(\alpha^{(a)})$
 \begin{equation}\label{eq:delta}
   \frac{(-1)^{\binom{p}{2}}}{n^p}
   \sum_I  (x_1\cdots x_p) (x_1+\cdots+x_p)^{|\alpha^{(a)}|}
   \det(x_j^{n-\alpha_i}) \det(x_i^{p-j})\,.
 \end{equation}

Let $(k_1,\ldots,k_p)=J(\alpha^{(a)})$.
If we write $a=pl+r$ with $0\leq r<p$, then this is the sequence
$$
   (ln+\alpha_{r+1},\ldots,ln+\alpha_p,
     (l{+}1)n+\alpha_1,\ldots,(l{+}1)n+\alpha_r)\,.
$$
The vector $(x_i^{n-k_j})$ is
\medskip

\noindent\qquad$((-1)^{(p+1)l}x_i^{n-\alpha_{r+1}},\ldots,
             (-1)^{(p+1)l}x_i^{n-\alpha_p},$ \smallskip

\hfill$(-1)^{(p+1)(l+1)}x_i^{n-\alpha_1},\ldots,
      (-1)^{(p+1)(l+1)}x_i^{n-\alpha_r})$\,.\qquad\medskip

\noindent
Thus we see that
$$ 
  \det(x_i^{n-k_j})\ =\ (-1)^{(p+1)a+r(p-r)}
            \det(x_i^{n-\alpha_j})\ =\ \det(x_i^{n-\alpha_j})\,,
$$
since, as in the calculation at the end of Section 2,
$r(p-r)\equiv pa+a$ modulo $2$.

Since $|\alpha^{(a)}|=\sum_jk_j-j$, we may substitute the last formula
into~(\ref{eq:delta}) and obtain
 \begin{eqnarray*}
  \delta(\alpha^{(a)})&=&
   \frac{(-1)^{\binom{p}{2}}}{n^p}\sum_I (x_1\cdots x_p)
      (x_1+\cdots+x_p)^{\sum_jk_j-j}
      \det(x_j^{n-k_i})\,\det(x_i^{p-j})\\
   &=& D(J(\alpha^{(a)}))\,,
 \end{eqnarray*}
as claimed.

We complete the proof of Lemma~\ref{lem:funD}(ii) and hence of
Theorem~\ref{thm:last} with the calculation below.

\begin{lemma}\label{ToDo}
$$
   \det\left(\frac{\partial c_i}{\partial x_j}\right) \ =\ 
     \prod_{i<j}(x_i-x_j)\ =\ \Delta\,.
$$
\end{lemma}

\begin{proof}
Let $F_p(x_1,\ldots,x_p)$ be this determinant.
Since
$$
    \frac{\partial c_i}{\partial x_j}\ =\ 
     c_{i-1}(x_1,\ldots,\widehat{x_j},\ldots,x_p)\,,
$$
where $\widehat{x_j}$ indicates that $x_j$ is omitted, 
we seek the determinant of the matrix whose $(i,j)$th entry is
$c_{i-1}(x_1,\ldots,\widehat{x_j},\ldots,x_p)$.
If we subtract the first column from each of the rest, we obtain a matrix in
block form
$$
  \left(\begin{array}{cc}1&0\\ *&A\end{array}\right)\,,
$$
where the entries of $A$ in position $(i,j)$
(note the shift from the original matrix) are
$$
  c_i(x_1,\ldots,\widehat{x_{j+1}},\ldots,x_p)-
  c_i(\widehat{x_1},\ldots,x_p)\ =\ 
  (x_1-x_{j+1})c_{i-1}(x_2,\ldots,\widehat{x_{j+1}},\ldots,x_p)\,.
$$
Dividing the common factors of $(x_1-x_{j+1})$ from the columns of $A$ gives
the matrix with entries
$c_{i-1}(x_2,\ldots,\widehat{x_{j+1}},\ldots,x_p)$,
and so we have the recursive formula
$$
  F_p(x_1,\ldots,x_p)\ =\ \prod_{j=2}^p(x_1-x_j)\cdot
           F_{p-1}(x_2,\ldots,x_p)\,.
$$
Since $F_1(x_p)=1$, this completes the Lemma.
\end{proof}

\subsection{The quantum cohomology ring of the Grassmannian}

We discuss an alternative view of the quantum
cohomology ring of the Grassmannian, mention how this ring arose in
representation theory, and give some open problems.

The presentation~(\ref{eq:qcoh-pres}) of 
$QH^*({\it Grass}(p,{\mathbb C}^{m+p}))$ is not what one ordinarily sees in
algebraic geometry, but rather an integral form with a parameter $q$
$$
  {\mathbb Z}[q][c_1,c_2,\ldots,c_p]/
   \langle h_{m+1},  \ldots,
            h_{m+p-1}, h_{m+p} + (-1)^pq \rangle\,.
$$
Then the genus zero Gromov-Witten invariant 
$\langle X_1,X_2,\ldots,X_l\rangle_0$ of cycles represented by classes
$\xi_1,\xi_2,\ldots,\xi_l\in{\mathbb Z}[c_1,c_2,\ldots,c_p]$
is the coefficient of $c_p^m\cdot q^d$ in the product
$\xi_1\cdot\xi_2\cdots\xi_l$.
Here $d$ is the degree of the curves this invariant enumerates.
In this presentation, the variable $q$ (which is the variable $\beta$
of Section~\ref{sec:VI}) keeps track of the degrees of curves.
This ring is graded, if $q$ has cohomological degree 
$2(m+p)$.

The cohomology of the Grassmannian has a basis of Schubert classes,
$\sigma_\alpha$, given by the Giambelli formula
(Jacobi-Trudi for combinatorists).
 \begin{equation}\label{eq:Giam}
  \sigma_\alpha\ =\ \det (h_{\alpha_i-j})_{1\leq i,j\leq p}\,.
 \end{equation}
The class $\sigma_\alpha$ is Poincar\'e dual to the Schubert variety
$\Omega_{\alpha^\vee}$.
The quantum cohomology ring of ${\it Grass}(p,{\mathbb C}^{m+p})$
may be viewed additively as polynomials in $q$ with coefficients in 
the cohomology of the Grassmannian, 
$H^*({\it Grass}(p,{\mathbb C}^{m+p})) [q]$, 
but with a deformed product, $*$, defined by
 \begin{equation}\label{eq:q-prod}
  \sigma_\alpha * \cdots * \sigma_\beta \ =\ \sum_{\gamma,d\geq 0}\ q^d \;
    \langle \Omega_{\alpha^\vee},\ldots,\Omega_{\beta^\vee},\,
      \Omega_{\gamma}\rangle_0^d\ \sigma_\gamma\,,
 \end{equation}
where $\langle \, \cdots\, \rangle_0^d$ is the genus 0 Gromov-Witten invariant for
degree $d$ curves.
(The degree $d$ is determined by the cohomological degrees of the
Schubert classes.)

Bertram~\cite{Be97} studied this ring, and showed that the Giambelli
formula~(\ref{eq:Giam}) remains valid with the quantum multiplication.
He also established a Pieri formula
$$
  \sigma_\alpha * h_a\ =\ \sum \sigma_\beta + q \sum \sigma_\gamma\,,
$$
the sum over all $\beta,\gamma$ with $|\beta|=|\alpha|+a$
and $|\gamma|=|\alpha|+a-m-p$, and satisfying
 \begin{eqnarray*}
  &\alpha_1\leq\beta_1<\alpha_2\leq\beta_2<\cdots<\alpha_p\leq\beta_p\,,&\\
  &\gamma_1\leq\alpha_1-1<\gamma_2\leq\alpha_2-1<\cdots<
     \gamma_p\leq\alpha_p-1\,.
 \end{eqnarray*}
Like the classical Giambelli and Pieri formulae~\cite{MR48:2152}, these
determine the ring structure of quantum cohomology with respect to
the basis of Schubert classes.

In particular, the structure constants $N^\gamma_{\alpha\,\beta}(m,p)$
defined by the formula
 \begin{equation}\label{eq:q-lr}
  \sigma_\alpha *\sigma _\beta \ =\ 
   \sum_{\gamma,d} q^d N^\gamma_{\alpha\,\beta}(m,p) \sigma_\gamma
 \end{equation}
are completely determined.
(Here, the summation is over $\gamma,d$ with $|\gamma|+d(m+p)=|\alpha|+|\beta|$.)
These are the analogs of the classical Littlewood-Richardson coefficients.
Like them, these numbers $N^\gamma_{\alpha\,\beta}(m,p)$ are non-negative.
Unlike the classical coefficients, there is as yet no quantum
Littlewood-Richardson formula for these constants which proves their
non-negativity.
These are certain three-point Gromov Witten invariants, as 
combining~(\ref{eq:q-prod}) and~(\ref{eq:q-lr}) shows
$$
  N^\gamma_{\alpha\,\beta}(m,p)\ =\ 
  \langle \Omega_{\alpha^\vee},\Omega_{\beta^\vee},\,
      \Omega_{\gamma}\rangle_0\,.
$$

These are known in the case of the Pieri formula and when $q=0$; for then they are the classical
Littlewood-Richardson coefficients.
The only case for which there is such a positive formula
is due to Tudose~\cite{Tu00}, when the
minimum of $m$ or $p$ is 2.
A formula for $N^\gamma_{\alpha\,\beta}(m,p)$ which involves signs
(like the formula~(\ref{formula}) for $d(q;m,p)$) was given by
Bertram, Ciocan-Fontanine, and Fulton~\cite{BCF}.
Interestingly, a similar formula was given previously in two
different contexts.

The Verlinde algebra is a variant of the representation ring 
of ${\mathfrak s}{\mathfrak l}_p$ where the usual product is replaced by the
fusion product, which is the tensor product of two
representations reduced at level $m$.
Witten~\cite{Wi95} explained the isomorphism between the Verlinde algebra
and the quantum cohomology ring of the Grassmannian, and this was rigorously
established by Agnihotri~\cite{Ag}.
This isomorphism is an analog of the relation between the cohomology rings of 
the Grassmann varieties ${\it Grass}(p,{\mathbb C}^{m+p})$, as $m$ varies, and 
the representation ring of  ${\mathfrak s}{\mathfrak l}_p$.
A formula similar to that of Bertram, Ciocan-Fontanine, and Fulton was given by 
Kac~\cite[Exercise 13.35]{Ka} and Walton~\cite{Wa90} in this context, where further details
may be found.

The cohomology ring of the Grassmannian is likewise isomorphic to an external  
representation ring of the symmetric groups (see~\cite{Fu97,Sa91}).
Similarly, there is a family of quotients of Hecke algebras at a
primitive $(m+p)$th root of unity whose external representation ring
is isomorphic to the quantum cohomology of the Grassmannian.
This was studied by Goodman and Wenzl~\cite{GW90}, and they also gave a
formula for  $N^\gamma_{\alpha\,\beta}(p,m)$ identical to that of 
Kac and Walton.

They also gave another presentation of the quantum cohomology ring
$$
  {\mathbb Z}[c_1,c_2,\ldots,c_p]/I_{m,p}\,,
$$
where $I_{m,p}$ is the ideal generated by
$$
  \{S_K\mid k_p-k_1=m+p\}\,,
$$
here, $K:0<k_1<k_2<\cdots<k_p$ is an increasing sequence of positive integers
and $S_K$ is defined by the Giambelli formula
$$
  S_K\ :=\ (h_{k_i-j})_{1\leq i,j\leq p}\,,
$$
where the polynomials $h_i$ are defined recursively by~(\ref{eq:Hs}).
Goodman and Wenzl proved that this quotient ring has an integral basis consisting of 
the classes
$$
  \{ S_I\,\mid\, I:0<i_1<i_2<\cdots<i_p < i_1+m+p\}\,.
$$
This is just the set of sequences $J(\alpha^{(a)})$ for all 
$\alpha^{(a)}\in{\mathcal C}_{m,p}:=\bigcup_b{\mathcal C}^b_{m,p}$,
where $J(\,\cdot\,)$ is the map~(\ref{eq:seq_def}) of Section~\ref{sec:Formulae}.
The relation between this basis and the basis $q^a\sigma_\alpha$ of
the conomology ring $H^*({\it Grass}(p,{\mathbb C}^{m+p}))[q]$ with the
variable $q$ adjoined  is just 
$$
  q^a\sigma_\alpha\ =\ S_{J(\alpha^{(a)})}\,.
$$
Bertram's Pieri formula has a nice expression in this basis:
$$
  h_a * S_I\ =\ \sum S_J\,,
$$
the sum over all sequences $J$ with $|J|=|I|+a$, where
$$
  i_1\leq j_1< i_2\leq j_2 < \cdots < i_p\leq j_p < i_1+m+p\,.
$$

We close with two additional problems concerning the quantum Littlewood-Richardson
coefficients.

Let $f^{\alpha^{(a)}}=\deg Z_{\alpha^{(a)}}$ be the number of saturated
chains in the poset ${\mathcal C}_{m,p}$ from the minimal element to 
$\alpha^{(a)}$.
More generally, given $\beta^{(b)},\alpha^{(a)}\in{\mathcal C}_{m,p}$, let
$f^{\alpha^{(a)}}_{\beta^{(b)}}$ be the number of saturated chains that begin at
$\beta^{(b)}$ and end at $\alpha^{(a)}$.
Write $\sigma_{\alpha^{(a)}}$ for the class $q^a\sigma_\alpha$
and then the expansion~(\ref{eq:q-lr}) becomes
$$
  \sigma_{\alpha^{(a)}} * \sigma_{\beta^{(b)}}\ =\ 
   \sum_{|\gamma^{(c)}|=|\alpha^{(a)}|+|\beta^{(b)}|} 
    N^{\gamma}_{\alpha\,\beta}(m,p)\ \sigma_{\gamma^{(c)}}\,.
$$
Iterating the Pieri formula with $a=1$ implies that
$$
  h_1^{*l}*\sigma_{\beta^{(b)}}\ =\ \sum_{|\gamma^{(c)}|=l+|\beta^{(b)}|} 
   f^{\gamma^{(c)}}_{\beta^{(b)}}\ \sigma_{\gamma^{(c)}}\,.
$$
Since $h_1^{*l}\ =\ \sum_{|\alpha^{(a)}|=l} f^{\alpha^{(a)}}\sigma_{\alpha^{(a)}}$, we 
may expand the left hand side to obtain
$$
  \sum_{|\alpha^{(a)}|=l} f^{\alpha^{(a)}}
     \sigma_{\alpha^{(a)}}*\sigma_{\beta^{(b)}}
  \ =\ 
  \sum_{|\alpha^{(a)}|=l}\ \   \sum_{|\gamma^{(c)}|}
    \ f^{\alpha^{(a)}} N^{\gamma}_{\alpha\,\beta}(m,p)\ \sigma_{\gamma^{(c)}} \,.
$$
Equating the coefficients of $\sigma_{\gamma^{(c)}}$, we obtain
\begin{equation}
  \sum_{|\alpha^{(a)}|=l} f^{\alpha^{(a)}}N^{\gamma}_{\alpha\,\beta}(m,p)
   \ =\ f^{\gamma^{(c)}}_{\beta^{(b)}}\,.
\end{equation}
The eventual combinatorial formula for the quantum Littlewood-Richardson
coefficients should also explain this identity.
That is, there should be some algorithm to convert a path in the
poset  ${\mathcal C}_{m,p}$ from $\beta^{(b)}$ to $\gamma^{(c)}$ into a path of the same 
length that starts at the minimal element, where the multiplicity of the occurrence of any
path to $\alpha^{(a)}$ is the quantum Littlewood-Richardson coefficient
$N^{\gamma}_{\alpha\,\beta}(m,p)$.
In short, we ask for a quantum version of Schensted insertion.

Another open problem is to show the (apparent) inequalities:
$$
  N^{\gamma}_{\alpha\,\beta}(m,p)\ \leq\ 
  N^{\gamma}_{\alpha\,\beta}(m+1,p)\ \leq\ 
  N^{\gamma}_{\alpha\,\beta}(m+2,p)\ \leq \ \cdots\ ,
$$
which were conjectured by Walton~\cite{Wa90}.

\section*{Acknowledgements}
We thank Joachim Rosenthal, who taught us the basics of systems theory and
commented on an early version of this manuscript, Jan
Verschelde for the matrix manipulations of Section 2.2, Eduardo Cattani
who resolved some of our questions on residues, and Emma Previato,
who solicited this survey.
We also thank Anders Buch, Sergey Fomin, Christian Lenart, Sasha Postnikov, Bruce Sagan,
Mark Shimozono, and Richard Stanley who each provided us with a last-minute proof of the 
identity~(\ref{eq:power_sum}). 

\providecommand{\bysame}{\leavevmode\hbox to3em{\hrulefill}\thinspace}
\providecommand{\MR}{\relax\ifhmode\unskip\space\fi MR }
\providecommand{\MRhref}[2]{%
  \href{http://www.ams.org/mathscinet-getitem?mr=#1}{#2}
}
\providecommand{\href}[2]{#2}

\end{document}